\documentclass[a4paper,10pt]{amsart}
\usepackage{graphicx}

\usepackage{amssymb}
\usepackage{amsthm}
\usepackage{latexsym}

\usepackage{color} 
\newcommand\ec { \color{black}}%
\newtheorem{thm}{Theorem}
\newtheorem{prop}{Proposition}
\newtheorem{lemma}{Lemma}
\newtheorem{cor}{Corollary}
\theoremstyle{definition}

\newtheorem{ex}{Example}
\newtheorem{exs}{Examples}
\newtheorem{rem}{Remark}
\newtheorem{rem-q}[thm]{Remark and Open Problems}
\newcommand{\fbar}{\boldsymbol{\overline{F}}}

\newcommand{\F}{\mathcal{F}}

\newcommand{\til}{{\widetilde{\star}}}
\newcommand{\M}{\mathcal{M}}

 \DeclareMathOperator{\Spec}   {Spec}
 
 \DeclareMathOperator{\Inv}{Inv}

 \newcommand{\ff}{\boldsymbol{f}}
  \newcommand{\fF}{\boldsymbol{F}}

\begin{document}

\title[ Pr\"ufer $\star$--multiplication domains and $\star$--coherence  ]{Pr\"ufer $\star$--multiplication domains and $\star$--coherence}
\author{Marco Fontana \hskip 1cm Giampaolo Picozza}
\address{Dipartimento di Matematica, Universit\`a degli Studi ``Roma Tre'',  \newline $\,$ \hskip 0.26 cm Largo San Leonardo Murialdo 1, 00146 Roma }

\email{fontana@mat.uniroma3.it, \  picozza@mat.uniroma3.it}

\thanks{This research was partially supported by the MIUR, under Grant PRIN 2005-015278.}


\date{\today}

\dedicatory{\rm Dipartimento di Matematica, Universit\`a degli Studi ``Roma Tre''}%

\maketitle

\section{ Introduction and Background}

The purpose of this paper is to deepen the study of the Pr\"ufer
$\star$--multiplication domains, where $\star$ is a semistar
operation (the definitions are recalled later in this section).
For this reason,  in Section 2,  we introduce the
$\star$--domains, as a natural extension of the $v$--domains
\cite[page 418]{gilmer}, where $v$ is the classical Artin's
divisorial operation.  We investigate their close relation with
the Pr\"ufer $\star$--multiplication domains. In particular, in
Section 3, we obtain a  characterization of Pr\"ufer
$\star$--multiplication domains in terms of $\star$--domains
satisfying a variety of  coherent-like conditions.  In Section 4,
we extend to the semistar setting the notion of
$\texttt{H}$--domain introduced by Glaz and Vasconcelos
\cite[Remark 2.2 (c)]{GV77}  and we show, among the other results
that, in the class of   the $\texttt{H}(\star)$--domains, the
Pr\"ufer $\star$--multiplication domains coincide with the
$\star$--domains.
\smallskip

Let $D$ be an integral domain with quotient field $K$. Let
$\boldsymbol{\overline{F}}(D)$ denote the set of all nonzero
$D$--submodules of $K$ and let $\boldsymbol{F}(D)$ be the set of
all nonzero fractional ideals of $D$, i.e. $E \in
\boldsymbol{F}(D)$ if $E \in \boldsymbol{\overline{F}}(D)$ and
there exists a nonzero $d \in D$ with $dE \subseteq D$. Let
$\boldsymbol{f}(D)$ be the set of all nonzero finitely generated
$D$--submodules of $K$. Then, obviously $\boldsymbol{f}(D)
\subseteq \boldsymbol{F}(D) \subseteq
\boldsymbol{\overline{F}}(D)$.

A \emph{semistar operation} on $D$ is a map $\star:
\boldsymbol{\overline{F}}(D) \to \boldsymbol{\overline{F}}(D), E
\mapsto E^\star$,  such that, for all $x \in K$, $x \neq 0$, and
for all $E,F \in \boldsymbol{\overline{F}}(D)$, the following
properties hold:
\begin{enumerate}
\item[$(\star_1)$] $(xE)^\star=xE^\star$;
 \item[$(\star_2)$] $E
\subseteq F$ implies $E^\star \subseteq F^\star$;
\item[$(\star_3)$] $E \subseteq E^\star$ and $E^{\star \star} :=
\left(E^\star \right)^\star=E^\star$.
\end{enumerate}

 Recall that,  given a
semistar operation $\star$ on $D$, for all $E,F \in
\boldsymbol{\overline{F}}(D)$,  the following basic formulas
follow easily from the axioms:    $$\begin{array}{rl} (EF)^\star
=& \hskip -7pt (E^\star F)^\star =\left(EF^\star\right)^\star
=\left(E^\star
F^\star\right)^\star\,;\\
(E+F)^\star =& \hskip -7pt \left(E^\star + F\right)^\star= \left(E
+
F^\star\right)^\star= \left(E^\star + F^\star\right)^\star\,;\\
(E:F)^\star \subseteq & \hskip -7pt (E^\star :F^\star) = (E^\star
:F) =
\left(E^\star :F\right)^\star,\;\, \mbox{\rm if \ } (E:F) \neq 0;\\
(E\cap F)^\star \subseteq & \hskip -7pt E^\star \cap F^\star =
\left(E^\star \cap F^\star \right)^\star,\;\, \mbox{\rm if \ }
E\cap F \neq (0)\,;
\end{array}
$$
 \noindent cf.  for instance \cite[Theorem 1.2 and p.  174]{FH2000}.

A \emph{(semi)star operation} is a semistar operation that,
restricted to $\boldsymbol{F}(D)$,  is a star operation (in the
sense of \cite[Section 32]{gilmer}). It is easy to see that a
semistar operation $\star$ on $D$ is a (semi)star operation if and
only if $D^\star = D$.

 If $\star$ is a semistar operation on $D$, then we can
consider a map\ $\star_{\!_f}: \boldsymbol{\overline{F}}(D) \to
\boldsymbol{\overline{F}}(D)$ defined, for each $E \in
\boldsymbol{\overline{F}}(D)$, as follows:

\centerline{$E^{\star_{\!_f}}:=\bigcup \{F^\star\mid \ F \in
\boldsymbol{f}(D) \mbox{ and } F \subseteq E\}$.}

\noindent It is easy to see that $\star_{\!_f}$ is a semistar
operation on $D$, called \emph{the semistar operation of finite
type associated to $\star$}.  Note that, for each $F \in
\boldsymbol{f}(D)$, $F^\star=F^{\star_{\!_f}}$.  A semistar
operation $\star$ is called a \emph{semistar operation of finite
type} if $\star=\star_{\!_f}$.  It is easy to see that
$(\star_{\!_f}\!)_{\!_f}=\star_{\!_f}$ (that is, $\star_{\!_f}$ is
of finite type).

If $T$ is an overring of $D$, we can define a semistar operation
on $D$, denoted by $\star_{\{T\}}$ and defined by
$E^{\star_{\{T\}}} := ET$, for each $E \in
\boldsymbol{\overline{F}}(D)$. It is easily seen that
$\star_{\{T\}}$ is  a semistar (non (semi)star,  if $D\subsetneq
T$) operation  of finite type.

If $\star_1$ and $\star_2$ are two semistar operations on $D$, we
say that $\star_1 \leq \star_2$ if $E^{\star_1} \subseteq
E^{\star_2}$, for each $E \in \fbar(D)$. This is equivalent to say
that $\left(E^{\star_{1}}\right)^{\star_{2}} = E^{\star_2}=
\left(E^{\star_{2}}\right)^{\star_{1}}$, for each $E \in
\fbar(D)$.  Obviously, for each semistar operation $\star$, we
have $\star_{\!_f} \leq \star$.

We say that a nonzero ideal $I$ of $D$ is a
\emph{quasi-$\star$-ideal} if $I^\star \cap D = I$, a
\emph{quasi-$\star$-prime} if it is a prime quasi-$\star$-ideal,
and a \emph{quasi-$\star$-maximal} if it is maximal in the set of
all quasi-$\star$-ideals. A quasi-$\star$-maximal ideal is  a
prime ideal. It is possible  to prove that each quasi-$\star_{_{\!
f}}$-ideal is contained in a quasi-$\star_{_{\! f}}$-maximal
ideal.  More details can be found in \cite[page 4781]{FL03}. We
will denote by $\M(\star_{_{\! f}})$ the set of the
quasi-$\star_{_{\! f}}$-maximal ideals of $D$.

  If $\Delta$ is a set of prime ideals of an integral
  domain
  $D$,  then the semistar operation $\star_\Delta$ defined
  on $D$ as follows

  \centerline{$
  E^{\star_\Delta} := \bigcap \{ED_P \;|\;\, P \in \Delta\}\,,
  \;  \textrm {  for each}    \; E \in \boldsymbol{\overline{F}}(D)\,,
  $}

  \noindent is called \it the spectral semistar operation associated to
  \rm
  $\Delta$.
  A semistar operation $\star$ of an integral domain $D$ is
  called
  \it a
  spectral semistar operation \rm if there exists a subset $
  \Delta$ of the prime spectrum of $D$, $\mbox{\rm Spec}(D)$,  such that $\,\star =
  \star_\Delta\,$.

  When $\Delta := \M(\star_{_{\! f}})$, we set $\til:= \star_{ \M(\star_{_{\! f}})}$, i.e.

  \centerline{$E^\til := \bigcap \left \{ED_M \, \vert \, M \in
\mathcal{M}(\star_{_{\! f}})\right\}$,  \;  for each $E \in \boldsymbol{\overline{F}}(D)$.}

A semistar operation $\star$ is \emph{stable} if $(E \cap F)^\star
= E^\star \cap F^\star$, for each $E,F \in \fbar(D)$.
 Spectral
semistar operations are stable
\cite[Lemma~4.1 (3)]{FH2000}.

We recall from \cite[Chapter~V]{FHP97}  that \it a localizing
system of ideals of $D$ \rm  is a family $\F$ of ideals of $D$
such that:

\begin{enumerate}
\item[{\bf(LS1)}] If  $I \in \F$ and $J$ is an ideal of $D$ such
that $I \subseteq J$, then $J \in \F$.
  \item[{\bf(LS2)}] If $I \in \F$ and $J$ is an ideal of $D$ such that
  $(J:_DiD) \in \F$, for each $i \in I$, then $J \in \F$.
\end{enumerate}

A localizing system $\F$ is \emph{finitely generated} if, for each
$I \in \F$, there exists a finitely generated ideal $J \in \F$
such that $J \subseteq I$.

The relation between  stable semistar operations and localizing
systems has been deeply investigated by M. Fontana and J. Huckaba
in \cite{FH2000} and by F. Halter-Koch in the context of module
systems  \cite{HK01}.   We summarize  some of results that we need
\rm in the following Proposition (see \cite[Proposition~2.8,
Proposition~3.2, Proposition~2.4, Corollary~2.11, Theorem~2.10
(B)]{FH2000}).

\begin{prop}
 \label{prop:loc1}  Let $D$ be an integral domain.
\begin{enumerate}
 \item[\bf (1)\rm] If $\star$ is a semistar operation on $D$, then  $\F^\star:=\left\{I \mbox{ ideal of $D$ } \mid I^\star =
D^\star \right\}$ is a localizing system (called \emph{the
localizing system associated to $\star$}).
 \item[\bf (2) \rm] If $\star$ is a semistar operation of finite type, then
$\F^\star$ is a finitely generated localizing system.
  \item[\bf (3) \rm] Let $\star_\F$ or, simply, ${\overline{\star}}$ be the semistar operation associated to a given  localizing system $\F$ of $D$ and defined by $E \mapsto E^{\overline{\star}}:=
\bigcup \left \{(E:J) \,\vert\, J\in\F \right \}$, for each $E \in
\fbar(D)$.  Then  $\star_\F$ (called \emph{the
semistar operation  associated to the localizing system $ \F$}) is  a stable semistar operation on $D$.
  \item[\bf (4) \rm] ${\overline{\star}} \leq \star$ and $\F^{\star} = \F^{\overline{\star}}$.
    \item[\bf (5) \rm] $\overline{\star} = \star$ if and only if $\star$ is
  stable.
  \item[\bf (6) \rm] If $\F$ is a finitely
generated localizing system, then $\star_\F$ is a finite type
(stable) semistar operation.

  \item[\bf (7) \rm]  $\F^{\star_{_{\! f}}}= (\F^\star){_{_{\! f}}} := \{ I \in \F^\star \mid I\supseteq J, \mbox{ for some finitely generated ideal } J\in \F^\star\} $ and  $\widetilde{\star}=\overline{\ \!\star_{_{\! f}}}$, i.e. $ \widetilde{\star}$ is the stable semistar operation of finite type associated to the localizing
system $\F^{\star_{_{\! f}}}$.
   \item[\bf (8) \rm] If $\F^\prime$ and $\F^{\prime \prime}$ are
two localizing systems of $D$, then $\F^\prime \subseteq
\F^{\prime \prime}$ if and only if $\star_{_{\F^\prime}} \leq
\star_{_{\F^{\prime \prime}}}$. \hfill $\Box$
\end{enumerate}
\end{prop}

 By $v_D$ (or, simply, by $v$) we denote  the $v$--(semi)star
operation defined as usual by  $E^v := (D:(D:E))$, for each $E\in
\boldsymbol{\overline{F}}(D)$. By  $t_D$ (or, simply, by $t$) we
denote  $(v_D)_{_{\! f}}$ the $t$--(semi)star operation on $D$ and
by  $w_D$ (or just by $w$) the stable semistar operation of finite
type associated to $v_D$ (or, equivalently, to $t_D$), considered
by Wang Fanggui and R.L. McCasland in \cite{WM97};  i.e. $w_D :=
\widetilde{v_D} = \widetilde{t_D}$.

 If $I \in \boldsymbol{\overline{F}}(D)$, we say that $I$ is
\emph{$\star$--finite} if there exists $J \in \boldsymbol{f}(D)$
such that $J^\star=I^\star$.  It is immediate to see that if
${\star_1} \leq {\star_2}$ are semistar operations and $I$ is
${\star_1}$--finite, then $I$ is ${\star_2}$--finite. In
particular, if $I$ is $\star_{\!_f}$--finite, then it is
$\star$--finite. The converse  is not true and it is possible to
prove that $I$ is $\star_{\!_f}$--finite if and only if there
exists $J \in \boldsymbol{f}(D)$, $J \subseteq I$, such that
$J^\star=I^\star$  \cite[Lemma 2.3]{FP}.

If $I$ is a nonzero ideal of $D$, we say that $I$ is
\emph{$\star$--invertible} if $(II^{-1})^\star = D^\star$, i.e.,
if $II^{-1} \in \F^\star$.  We denote by $\Inv(D, \star)$ the set
of all the $\star$--invertible ideals of $D$.   From the
definitions, it follows easily   that an ideal is
$\star$--invertible if and only if it is
$\overline{\star}$--invertible (and so $I$ is
$\tilde{\star}$--invertible if and only if $I$ is $\star_{_{\!
f}}$--invertible).  If $I$ is $\star_{_{\! f}}$--invertible, then
$I$ and $I^{-1}$ are $\star_{_{\! f}}$--finite \cite[Proposition
2.6]{FP}.

A domain $D$ is called \it a  Pr\"ufer $\star$--multiplication
domain \rm (for short, \it P$\star$MD\rm) if each nonzero finitely
generated ideal is $\star_{_{\! f}}$-invertible  (cf. for instance
\cite{HMM} and \cite{FJS03b}).  When $\star=v$ we have the
classical notion of P$v$MD   (cf. for instance \cite{Gr},
\cite{MZ} and \cite{K89b});   when $\star =d$, where $d$ denotes
\it the identity (semi)star operation, \rm we have the notion of
Pr\"ufer domain   \cite[Theorem 22.1]{gilmer}.

 We say that   a semistar operation $\star$  on $D$   is \emph{a.b.} (= {\sl{arithmetisch brauchbar}}) if, for each $E \in
\boldsymbol{f}(D)$ and for all $F,G \in \fbar(D)$, $(EF)^\star
\subseteq (EG)^\star$ implies $F^\star  \subseteq  G^\star$ and
$\star$ is \emph{e.a.b.} (= {\sl{endlich arithmetisch brauchbar}})
if the same holds for all $E, F, G \in \boldsymbol{f}(D)$.
Obviously, a.b. implies e.a.b.; in case of semistar operations of
finite type, it is easy to see that the notions of e.a.b. and a.b.
semistar operation coincide (in this situation,  we will write
(e.)a.b. operation).

 Finally, let $\star_1 \leq \star_2$ be two semistar operations on
$D$, then we say that $D$ is \it  a $(\star_1, \star_2)$--domain \rm  if $\star_1
= \star_2$.

\bigskip

\bigskip

\section{Pr\"ufer $\star$--multiplication domains and $\star$--domains}
Let $D$ be a domain and $\star$ a semistar operation on $D$. We
say that $D$ is \it a $\star$--domain \rm if $(II^{-1})^\star = D^\star$
for each $I \in \ff(D)$.

For $\star = v$ we have the notion of $v$--domain considered in
\cite[Section 34]{gilmer}, and for $\star=d$  we have that the
notions of $d$--domain, P$d$MD and Pr\"ufer domain coincide
\cite[Theorem 22.1]{gilmer}.

\begin{prop}  \label{basic}
 Let $D$ be a domain and let $\star_1,  \star_2$ be two
semistar operations on $D$.
\begin{enumerate}
\item[\bf (1) \rm]  If  $\star_1 \leq \star_2$ and $D$ is a $\star_1$--domain then $D$
is a $\star_2$--domain.
\item[\bf (2) \rm] The following statements are equivalent:
\begin{enumerate}
    \item[\bf (i) \rm] $D$ is a $\star_{_{\! f}}$--domain;
   \item[\bf (ii) \rm] $D$ is a $\widetilde{\star}$--domain;
\item[\bf (iii) \rm] $D$ is a P$\star$MD.
\end{enumerate}
    \item[\bf (3) \rm] A P$\star$MD is always a  $\star$--domain.
      \item[\bf (4) \rm]  The following statements are equivalent:
\begin{enumerate}
  \item[\bf (j) \rm] $D$ is a $\star$--domain;
   \item[\bf (jj) \rm] $D$ is a $\overline{\star}$--domain.
\end{enumerate}
\item[\bf (5) \rm]  Let $v(D^\star)$ be  the semistar operation on
$D$ defined by $E \mapsto E^{v(D^\star)}:= (D^\star :
(D^\star:E))$, for each $E \in \fbar(D)$  (cf. \cite[Lemma 2.11
(4) and its proof]{FP} or \cite[Example 1.8 (2)]{GP}).  If $D$ is
a $\star$--domain, then $D$ is a $v(D^\star)$--domain and
$\star_{_{\! f}} = (v(D^\star))_{_{\! f}}$.
\end{enumerate}
\end{prop}
\begin{proof} (1) follows immediately from the definitions.

(2)  follows from \cite[Proposition 2.18 ]{FP} and from the definition of a P$\star$MD. (Note that a  P$\star$MD coincides with a  P$\widetilde{\star}$MD and with a P$\star_{_{\! f}}$MD, cf. also \cite[page 30]{FJS03b}.)

 (3)  is an easy consequence of (1) and (2),  since $\star_{_{\! f}} \leq \star$.

(4) We have already observed that $\mathcal F^\star = \mathcal
F^{\overline{\star}}$, thus $II^{-1} \in \mathcal F^\star $ if and
only if    $II^{-1} \in \mathcal F^{\overline{\star}} $.

 (5) Since $\star \leq v(D^\star)$,  for each semistar operation
$\star$ \cite[Corollary 3.8]{GP},   the first assertion is
an immediate consequence of (1). The second assertion follows by
\cite[Remark 2.13 (c)]{FP}.  \end{proof}

  It is known that, when $\star$ is a star operation, a
P$\star$MD is a P$v$MD such that $\star_{_{\! f}} =t$
\cite[Proposition 3.4]{FJS03b}. The next result extends the
previous characterization to the case of $\star$--domains.

 \begin{cor} Let $\star$ be a star operation defined on an
integral domain $D$. If $D$ is a $\star$--domain then $D$ is a
$v$--domain and $\star_{_{\! f}} = t$.
\end{cor}
\begin{proof}  Since in the present situation $D^\star =D$, the statement is a straightforward consequence of Proposition \ref{basic} (5).
\end{proof}

\begin{rem} \label{rk:2.2} \bf (1) \rm As a consequence of the previous result we re-obtain that the notions of P$v$MD, P$t$MD and P$w$MD coincide (cf. \cite[Theorem 2.18]{AC2000} and \cite[Remark 3.1 and Corollary 3.1]{FJS03b}).

 \bf (2) \rm Note that  $\star$--domains are not always P$\star$MD, even if $\star$ is a (semi)star operation.  For instance,  recall that an essential domain is a $v$--domain \cite[Proposition 44.13]{gilmer} and not every essential domain is a P$v$MD \cite{HO}  (cf. also \cite{H81}  for an example of an essential domain with a non-essential localization, and so, in particular,  which is not P$v$MD \cite[Example 2.1, Proposition 1.1 and Corollary 1.4]{MZ}).     An example of a star operation  $\star$, defined   on an essential domain $D$,  such that $D$ is a $\star$--domain but not a  P$\star$MD is given in the following Example \ref{b-domain}.

 \bf (3) \rm Note also that, from Propostion \ref{basic} (2) and the previous observation (2), we deduce in particular that the notions of $t$--domain and $w$--domain coincide,  but they are strictly stronger than the notion of $v$--domain (as observed in (2)).

 \bf (4) \rm   Note that from Proposition \ref{basic} (1, 2), we deduce that if   $\star_1 \leq   \star_2$ and $D$ is a P$\star_1$MD then  $D$ is also a P$\star_2$MD.  Since $\til \leq \overline{\star} \leq \star$,  we have that a P$\til$MD is a P$\overline{\star}$MD, which is  a P$\star$MD and thus  it is easy to see that all these notions coincide (cf. Proposition \ref{basic} (2)).

 \bf(5) \rm
In \cite[Section 34]{gilmer}, a $v$--domain is defined as a domain
such that the $v$--operation is e.a.b., and in \cite[Theorem
34.6]{gilmer} it is shown that this is equivalent to the fact that
each finitely generated ideal is $v$--invertible. This  type of
 characterization does not hold for general semistar operations
$\star$ (cf. the following Example \ref{eab-not-star}).
\end{rem}

\smallskip

\begin{prop}  \label{ab} Let $D$ be an integral domain and $\star$ a semistar operation on $D$. If $D$ is a $\star$--domain (e.g. a P$\star$MD)  then $\star$ is an a.b. operation. In particular, a $\star$--domain $D$ is  quasi--$\star$--integrally closed (i.e. $D^\star = \bigcup \{ (F^{\star} : F^{\star}) \mid F \in \ff(D) \}$) and so $D^\star $ is integrally closed.
\end{prop}
\begin{proof} If $(FG)^\star \subseteq (FH)^\star $ with $F \in \ff(D)$ and $G, H \in \fbar(D)$, then $G^\star = ((FF^{-1})^\star G)^\star $ $= (FF^{-1}G)^\star \subseteq  (FF^{-1}H)^\star= ((FF^{-1})^\star H)^\star= H^\star $.  The other statements follow respectively by \cite[Lemma 4.13]{EFP} and by    \cite[Proposition 4.3]{FL01a}.\end{proof}

Note that it is not always true that if $\star$ is a.b. then $D$
is a $\star$--domain, as the following examples show.

\begin{ex} \label{eab-not-star} \bf (1)  \sl An a.b.  non-stable semistar (non-star) operation of finite type $\star$ on an integral  domain $D$ such that $D$ is not  a $\star$--domain (or, equivalently, not a P$\star$MD). \rm

Take a pseudo-valuation non-valuation domain
$D$, with maximal ideal $M$, and set $V:=M^{-1}$. Let $\star =
\star_{f} := \star_{\{V\}}$. It is easy to see that $\star$ is an
 a.b. semistar operation on $D$.  Since $\M(\star_{f} ) = \{M\}$, then $D$ is not a P$\star$MD \cite[Theorem 1.1 ((1)$\Leftrightarrow$(4))]{HMM} and so it is not a
$\star$--domain, because in this case $\star$ is of finite type, so $D$ is a
$\star$-domain if and only if it is a P$\star$MD (cf. Proposition \ref{basic} (2))). Finally $\star$ is not stable, because otherwise $\star = \til$ \cite[Corollary 3.9 (2)]{FH2000} and hence the fact that $\til$ is  (e.)a.b. implies that $D$ is a P$\star$MD \cite[Theorem 3.1]{FJS03b}.

\bf (2)  \sl An a.b.  non-stable star operation of finite type $\star$ on an integral  domain $D$ such that $D$ is not  a $\star$--domain. \rm

It is easy to check that the $v$--operation is e.a.b. if and only
if the $t$--operation is a.b.  (cf. also \cite[Definition 3.6 and
Lemma 3.9 (2)]{FL06}).   Therefore, a $v$--domain is an integral
domain such that the $t$--operation is a.b. (Remark \ref{rk:2.2}
(5)), but we already observed (Remark \ref{rk:2.2} (2)) that a
$v$--domain is not necessarily a $t$--domain (that is, a P$v$MD).
The non-stability of the operation $t$,  on a $v$--domain which is
not a P$v$MD, follows from the same argument as  in  the previous
example.
\end{ex}

\begin{ex} \label{b-domain} \sl An (a.b.)  star operation $\star$ on an integral domain $D$ such that $D$ is a $\star$--domain but not a P$\star$MD. \rm

Let $D$ be a domain and  let  $\{ V_\alpha \}$ be a nonempty set
of valuation overrings of $D$ which are essential for $D$ (that
is, $V_\alpha$ is the localization of $D$ at its center $P_\alpha$
in $D$). Consider the semistar operation $\star$ induced by this
set overrings ( i.e. $E^\star :=   \bigcap_\alpha EV_\alpha$, for
each $E \in \fbar(D)$;  thus $\star$ is a semistar non-(semi)star
operation on $D$ if $D \subsetneq \bigcap_{\alpha} V_\alpha$).
Let $I$ be a nonzero finitely generated ideal of $D$, then
$(II^{-1})^\star = \bigcap_\alpha (I(D:I))V_\alpha  =
\bigcap_\alpha I(D:I)D_{P_\alpha} = \bigcap_\alpha
(ID_{P_\alpha}(D_{P_\alpha}:ID_{P_\alpha})) = \bigcap_\alpha
IV_\alpha(V_\alpha : IV_\alpha) = \bigcap_\alpha V_\alpha =
D^\star$. Thus, each nonzero finitely generated ideal of $D$ is
$\star$--invertible and so $D$ is a $\star$--domain.

 Note that a similar argument shows that $\star$ is stable:
$(E\cap F)^{\star} =  \bigcap_\alpha (E\cap F)D_{P_\alpha} =
\bigcap_\alpha (ED_{P_\alpha}  \cap FD_{P_\alpha}) =
(\bigcap_\alpha ED_{P_\alpha})  \cap (\bigcap_\alpha
FD_{P_\alpha})  = E^\star \cap F^\star$, for all $E,F \in
\fbar(D)$
 (i.e. $\star = \overline{\star}$).

 Note also that a semistar operation defined by a family of valuation overrings (like the $\star$ defined above) is necessarily a.b..

Assume from now that $D=\bigcap_{\alpha} V_\alpha$. Note that, in
this case,   $\star$, defined on  $\fF(D)$,   is a star operation
on $D$, thus $\star \leq v$ and so $D$ is also a $v$--domain. By
Proposition \ref{basic}(5), we can deduce that $\star_{_{\! f}} =
t$. So, $D$ is a P$\star$MD (= P$\star_{_{\! f}}$MD)  if and only if it is
a P$v$MD (= P$t$MD).  We can conclude that if you choose $D$  not
to be a P$v$MD (such example exists (Remark \ref{rk:2.2} (2)),
then $D$ is a $\star$-domain (and a $v$--domain) which is not a
P$\star$MD (nor a a P$v$MD).

In this situation,  $\star$ may be not of finite type, since if
$\star = \star_{_{\! f}}$, then $\star =\overline{\star} =\til$,
thus  $\til$ would be a.b. and so $D$ would be a P$\star$MD
\cite[Theorem 3.1 ((v)$\Rightarrow$(i))]{FJS03b}.    Finally, note
that if $\overline{v} \neq v$, then necessarily $\star \lneq v$.

\end{ex}

\begin{prop} \label{E-inv}  Let $D$ be an integral domain and $\star$ a semistar operation on $D$.  The following are equivalent:
\begin{enumerate}
\item[\bf (i) \rm]  $D$ is a $\star$--domain [respectively: a P$\star$MD].

\item[\bf (ii) \rm]   for all  $E, F \in \ff(D)$  there exists $H \subseteq (E:F)$, $H \in \fF(D)$  [respectively: $H \in \ff(D)$],  such that $E^\star = (FH)^\star$.

\item[\bf (iii) \rm]   for all  $E, F \in \ff(D)$, $(F(E:F))^\star = E^\star$  [respectively: $(F(E:F))^{\star _f}= E^\star$].
\end{enumerate}
\end{prop}
\begin{proof}  (i)$\Rightarrow$(ii) Assume that $D$ is a $\star$--domain and take $H:= F^{-1}E$. Clearly $HF = FF^{-1}E \subseteq DE = E$ and so $H\subseteq (E:F)$.   Moreover, $(FH)^\star = (FF^{-1}E)^\star = ((FF^{-1})^\star E)^\star = E^\star$.   If $D$ is a P$\star$MD, let $G\in \ff(D)$ such that $G \subseteq F^{-1}$ and $G^\star = (F^{-1})^\star$. In this case, we just need to modify the choice of $H$, setting $H:= GE$.

 (ii)$\Rightarrow$(iii) is straightforward, since $H\subseteq (E:F)$ [and, in the parenthetical statement, $H \in \ff(D)$].

  (ii)$\Rightarrow$(iii) is obvious by taking $E=D$.
  \end{proof}

  \begin{rem}  Note that the proof of Proposition \ref{E-inv} shows that  \it $D$ is a $\star$--domain [respectively: P$\star$MD] if and only if $(F(E:F))^\star = E^\star$ [respectively:  $(F(E:F))^{\star_{_{\! f}}} = E^{\star_{_{\! f}}}$], for all  $ F \in \ff(D)$ and $E \in \fbar(D)$.\rm
   \end{rem}

   We are in condition to give a characterization of the $\star$--domains [respectively: Pr\"ufer $\star$--multiplication domains] by using that $\star$ is a.b. or that $D$ is quasi--$\star$--integrally closed (Proposition \ref{ab}).

   \begin{cor} \label{cor:2.8} Let $D$ be an integral domain and $\star$ a semistar operation on $D$.  The following are equivalent:
\begin{enumerate}
\item[\bf (i) \rm]  $D$ is a $\star$--domain [respectively: a P$\star$MD].

 \item[\bf (ii) \rm]     $\star$ is a.b. and $(EF^{-1})^\star = (E^\star :F)$ [respectively:
 $(EF^{-1})^{\star_{_{\! f}}} = (E^{\star_{_{\! f}}} :F)$] for all $ F \in \ff(D)$ and $E \in \fbar(D)$.

  \item[\bf (iii) \rm]     $D$ is a quasi--$\star$--integrally closed domain and $(EF^{-1})^\star = (E^\star :F)$ [respectively:
 $(EF^{-1})^{\star_{_{\! f}}} = (E^{\star_{_{\! f}}} :F)$] for all $ F \in \ff(D)$ and $E \in \fbar(D)$.

  \end{enumerate}

   \end{cor}
   \begin{proof}  We show the equivalences for the $\star$--domain case; the equivalences among the parenthetical statements follow from the fact that a P$\star$MD coincide with a $\star_{_{\! f}}$--domain (Proposition \ref{basic} (2)).

   (i)$\Rightarrow$(ii)$\Rightarrow$(iii).  By Proposition \ref{ab} and by the fact that $E^\star = (FF^{-1}E)^\star \subseteq  (F(E:F))^\star   \subseteq  (F(E^\star:F))^\star \subseteq  E^\star$ we have $(FF^{-1}E)^\star = (F(E^\star:F))^\star$. Since $\star$ is a.b.  we deduce that $(F^{-1}E)^\star = (E^\star:F)^\star = (E^\star:F)$.

    (iii)$\Rightarrow$(i) By taking $E=F$ we have that $(FF^{-1})^\star =(F^\star :F)$. Moreover, by the fact that  $D$ is  quasi--$\star$--integrally closed,  we have $D^\star \subseteq (F^\star :F)= (F^\star :F^\star) \subseteq \bigcup \{ (F^\star :F^\star) \mid F\in \ff(D)\} = D^\star$, and so $(F^\star :F) = D^\star$.
    \end{proof}

    \smallskip

    Note that the previous corollary generalizes to the semistar setting some cha\-rac\-terizations of the $v$--domains proved in \cite[Theorem 2]{AACDMZ}.

    \medskip

    The next goal is to relate $\star$--domains with properties of stability for the semistar operation $\star$.

      \smallskip

    \begin{prop} \label{finite-stable} Let $D$ be a $\star$--domain.  Assume that $D$ is integrally closed. Then, for all $E,F \in \ff(D)$:
    $$ (E:_DF)^\star = (E^\star :_{D^\star} F)\,.$$
    In particular, $D$ is a $((\overline{\star})_{_{\! f}}, \star_{_{\! f}})$--domain and so $(E \cap F)^\star = E^\star \cap F^\star$.
    \end{prop}
    \begin{proof}   If $D$ is integrally closed, then for each $E \in \ff(D)$, $(E:E) =D$    \cite[Proposition 34.7]{gilmer}.
    Note that $  (E:_DF) =(E:F) \cap D = (E:F) \cap (E:E) = (E : F+E)$.  On the other hand,  we know already that in a $\star$--domain,   $(E^\star : E^\star) =D^\star$ and  $(EF^{-1})^\star =(E^\star : F)$ (Proposition \ref{ab} and Corollary \ref{cor:2.8}).  Since, in general, $(EF^{-1})^\star \subseteq (E :F)^\star \subseteq  (E^\star : F)$, we deduce that  $(E :F)^\star =  (E^\star : F)$.  Moreover
    $(E : F +E)^\star = (E^\star : F+E )=(E^\star : (F+E)^\star) = (E^\star : (F^\star+E^\star)^\star) = (E^\star : (F^\star+E^\star)) = (E^\star : F^\star)\cap (E^\star : E^\star)= (E^\star : F^\star)\cap D^\star=  (E^\star :_{D^\star} F^\star) = (E^\star :_{D^\star} F)$.

    In order to prove the second part of the statement we show that, for each $E \in \ff(D)$, $E^\star = E^{\overline{\star}} $.
    $$
    \begin{array}{rl}
     x \in E^\star \Leftrightarrow&
    1\in (E^\star: xD)  \Leftrightarrow  (E^\star:_{D^\star} xD) = D^\star \Leftrightarrow  (E:_{D} xD)^\star= D^\star \\  \Leftrightarrow & (E:_{D} xD)^\star \in \F^\star     \Leftrightarrow   I \subseteq (E: xD)   \,\; \mbox{ for some } I \in  \F^\star \\ \Leftrightarrow&
xI \subseteq E\,\; \mbox{ for some } I \in  \F^\star   \Leftrightarrow x \in E_{\F^\star} = E^{\overline{\star}}\,.
\end{array}$$
Finally, since ${\overline{\star}}$ is stable  and $
(\overline{\star})_{_{\! f}}= \star_{_{\! f}}$,  then $E^\star
\cap F^\star= E^{\overline{\star}}\cap F^{\overline{\star}} =(E
\cap F)^{\overline{\star}} \subseteq (E \cap F)^\star \subseteq
E^\star \cap F^\star$, when $E,F \in \ff(D))$.
\end{proof}

\begin{rem}  \label {rk:2.10} \bf (1) \rm Note that, if $\star$ is a (semi)star operation on $D$ and $D$ is a $\star$--domain, then $D$ is integrally closed, thus the previous Proposition \ref{finite-stable}  generalizes to the semistar setting a result proved recently  by Anderson and Clarke \cite[Theorem 2.8]{AC2005}.

\bf (2) \rm  In relation with Proposition \ref{finite-stable}, we
remark that  it is possible to generalize in the semistar setting
a result proved by Anderson and Cook \cite[Theorem 2.6]{AC2000}.
 More precisely,  \it  if $\star$ is a semistar operation on an
integral domain, then the following conditions are equivalent:
\begin{enumerate}

\item[\bf (i) \rm]  $D$ is a $((\overline{\star})_{_{\! f}},
\star_{_{\! f}})$--domain  [respectively:  $(\overline{\star},
\star)$--domain].

\item[\bf (ii) \rm]   For all $E, F \in \ff(D)$   [respectively:
$E, F \in \fbar(D)$],  $(E \cap F)^\star = E^\star \cap F^\star$.

\item[\bf (iii) \rm]   For all $E, F \in \ff(D)$  [respectively:
$E \in \fbar(D)$, $F \in \ff(D)$],  $(E:_DF)^\star = (E^\star
:_{D^\star} F)$.

\item[\bf (iv) \rm]   For each $E\in \ff(D)$   [respectively: $E
\in \fbar(D)$],  and for each nonzero element $x\in K$, $(E:_D
xD)^\star = (E^\star :_{D^\star} xD)$.
\end{enumerate}
Clearly, if $\star$ is stable  (i.e. $\overline{\star} = \star$),
  then all the previous statements hold.  \rm

The implications (iii)$\Rightarrow$(iv)$\Rightarrow$(i)$\Rightarrow$(ii) are  essentially proved in Proposition \ref{finite-stable}.

(ii)$\Rightarrow$(iii).  If $F =x_1D+x_2D+\cdots +x_nD$, then
 $(E :_DF)^\star =
(E:_D (x_1D+x_2D+\cdots +x_nD))^\star = (\bigcap \{  E:_Dx_iD)\mid 1 \leq i\leq n\})^\star =
 (\bigcap \{  x_i^{-1}E \cap D\mid 1 \leq i\leq n\})^\star =  \bigcap \{  x_i^{-1}E^\star \cap D^\star\mid 1 \leq i\leq n\}= \bigcap \{  E^\star:_{D^\star} x_iD)\mid 1 \leq i\leq n\} = (E^\star :_{D^\star}F)$.

\bf (3) \rm It is easily seen that, \it  if $\star$ is a semistar
operation of finite type on $D$, $D$ is a
$((\overline{\star})_{_{\! f}}, \star_{_{\! f}})$--domain if and
only if $\star$ is a stable semistar operation. \rm Thus a
semistar operation of finite type is stable if and only if the
(non-parenthetical)  equivalent conditions in (2) are satisfied.

  \bf (4) \rm    As already observed in the star setting, if $D$ is an
integrally closed $((\overline{\star})_{_{\! f}}, \star_{_{\! f}})$--domain then
$D$ is not necessarily a $\star$--domain. (Take $D$ to be an
integrally closed non-Pr\"ufer domain and $\star = d$.)

However,  in the particular case that $\star = v$, then  \it  $D$
is a $v$--domain if and only if $D$ is integrally closed
$((\overline{v})_{_{\! f}},  t)$--domain. \rm  (The ``if part'' is
 due to Anderson et al.  \cite[Theorem 7]{AACDMZ}, cf. also
(1), (2)   and Proposition \ref{finite-stable}; the ``only if ''
part  was proved by Matsuda and Okabe \cite{MO1993},  cf. also
\cite[Theorem 2.8]{AC2005}.)

 At this point, for the general case, if we replace the condition
``$D$ is integrally closed'' with the condition ``$\star$ is a.b.
on $D$'' (which is a stronger condition in the (semi)star
setting), it is natural to ask:

\bf (Q-1) \rm
 \sl  Let $\star$ be a semistar operation on $D$. Is it true that $D$ is a $\star$--domain
if and only if $\star$ is a.b. and $D$ is a   $((\overline{\star})_{_{\! f}}, \star_{_{\! f}})$--domain ? \rm

 Note that  the  answer to the previous question is positive  for
$\star$ of finite type, since in this case $
(\overline{\star})_{_{\! f}} = \til$ and  we know that $D$ is a
P$\star$MD if and only if $\til = \star_{_{\! f}}$ is (e.)a.b.
\cite[Theorem 3.1]{FJS03b},  cf. also the following Theorem
\ref{star-domain-tilde}.

\ There is another important case in which  the answer to \bf
(Q-1) \rm is positive. Let $\star:= \star_\Delta$ be a spectral
semistar operation, where $\Delta \subseteq \Spec(D)$. Clearly
$\star$ is stable and so $D$ is a $((\overline{\star})_{_{\! f}},
\star_{_{\! f}})$--domain.

Assume that $\star$ is  a.b..  For each $P \in \Delta$, let
$\iota_{_{\! P}}$ be the canonical embedding of $D$ in $D_P$. We
claim that $\star_{\iota_P}$ coincides with $d_{D_P}$ (i.e. the
identity (semi)star operation of $D_P$), for each $P \in \Delta$.
In fact, if $E \in \fbar(D_P)$, $E \subseteq E^{\star_{\iota_{_{\!
P}}}} = E^\star = \bigcap_{P_\alpha \in \Delta} ED_{P_\alpha}
\subseteq ED_P = E$. Moreover, $\star_{\iota_{_{\! P}}} \
(=d_{D_P})$ is  also  a.b. by \cite[Proposition 3.1 (4)]{GP}.
Thus, each finitely generated ideal of $D_P$ is a cancellation
ideal and, so, $D_P$ is a valuation domain \cite[Theorem
24.3]{gilmer}. Therefore the semistar operation $\star$ is defined
by a family of valuation overrings of $D$ which are essential for
$D$. We have already shown in Example \ref{b-domain} that, in this
case, $D$ is a $\star$--domain.

Conversely,  if $D$ is a $\star$--domain then $\star$ is a.b.
(Proposition \ref{ab}) and, as we already remarked,  if $\star$ is
a spectral semistar operation then $\star$ is stable.

\end{rem}

  \bigskip

The next proposition generalizes to the case of $\star$--domains
some results already known for Pr\"ufer $\star$--multiplication
domains  (cf. \cite[Proposition 3.1 and 3.2]{FJS03b}).

\begin{prop} \label{prop:2.7} Let $T$ be an overring of an integral domain $D$ and let $\iota: D \hookrightarrow T$ be the canonical embedding.
\begin{enumerate}
\item[\bf (1) \rm]  Let $\star$ be a semistar operation on $D$ and let $\star_\iota$ be the semistar operation on $T$ defined by $E^{\star_\iota} := E^\star$,  for each $E \in \fbar(T) \ (\subseteq \fbar(D))$.  If $D$ is a $\star$--domain then $T$ is a $\star_\iota$--domain.
\item[\bf (2) \rm]  Let $\ast$ be a semistar operation on $T$ and let $\ast^\iota$ be the semistar operation on $D$ defined by $E^{\ast^\iota} := (ET)^\ast$,  for each $E \in \fbar(D)$.  If $T$ is a $\ast$--domain  and $\iota$ is flat then $D$ is a $\ast^\iota$--domain.
\end{enumerate}
\end{prop}
\begin{proof} \bf (1) \rm
 Let $G:= x_1T +x_2T+ \dots +x_n T \in \ff(T)$ and set $G_0:= x_1D+x_2D+ \dots +x_n D \ ( \in \ff(D)).$  Then $(G(T:G))^{\star_\iota} = (G_0T(T: G_0T))^\star \supseteq ((G_0 (D:G_0))T )^\star =
((G_0 (D:G_0))^\star T )^\star = (D^\star T)^\star  = T^\star$, thus  we conclude immediately  that  $(G(T:G))^{\star_\iota} $ coincides with $T^{\star_\iota}$.

\bf (2) \rm  Let $F:= x_1D+x_2D+ \dots +x_n D \in \ff(D)$. Then $(F(D:F))^{\ast^\iota} = ((F(D:F))T)^{\ast}  = (FT(T:FT))^{\ast}  = T^\ast = D^{\ast^\iota}$.
\end{proof}

 \smallskip

  \begin{rem}  Note that the semistar operation $v(D^\star)$ considered in Proposition \ref{basic}  (5) coincides with $(v_{D^\star})^\iota$ (notation as in Proposition \ref{prop:2.7} (2)), where $\iota: D \hookrightarrow D^\star$ is the canonical embedding and $v_{D^\star}$ is the $v$--(semi)star operation on $D^\star$.

  \end{rem}

\begin{ex}  \sl The assumption of flatness  is essential in  the proof of  Proposition \ref{prop:2.7} (2). \rm

Let $k \subset K$ be a proper finite extension of fields  and $X$
an indeterminate over $K$. Set $T:= K[X]_{(X)}$, $D : = k +
XK[X]_{(X)}$, $M:= XK[X]_{(X)}$,  $\iota: D \hookrightarrow T$ the
canonical embedding (which is clearly non-flat). Note that  $T$,
being a discrete valuation domain,  is a P$\ast$MD (and so a
$\ast$--domain) for all the semistar operations $\ast$ on $T$, in
particular $T$ is a P$d_T$MD, where $d_T$ is the identity
(semi)star operation on $T$. On the  other hand $D$ is not a
$(d_T)^\iota$--domain, since $(D:M) = (M:M) = T$, hence
$(MM^{-1})^{(d_T)^\iota} = (MT)^{(d_T)^\iota} = M \neq D$ and $M$
is finitely generated in $D$, by the finiteness of $\iota$
\cite[Proposition 1.8]{F}.

 As we have already observed (Example \ref{eab-not-star} (1)),
  $(d_T)^\iota$ is an a.b.  semistar  operation on $D$,
since $(d_T)^\iota =\star_{\{T\}}$ and $T$ is a valuation domain;
  therefore   $\star :=(d_T)^\iota \ (= \star_{\{T\}})$ gives an
example of an a.b. semistar operation on $D$ such that $D$ is not
a $\star$--domain.  Moreover,  $ d_D   = \widetilde{\
\star_{\{T\}}} = (\overline{\ \!\star_{\{T\}}})_f \lneq
(\star_{\{T\}})_f= \star_{\{T\}}$, since  $\star_{\{T\}}$  is
stable if and only if $\iota: D \hookrightarrow T$ is flat
(cf. \cite[Proposition 1.7]{Uda} and \cite[Theorem
7.4 (i)]{matsumura}).

\end{ex}

\smallskip The next result shows that a $\star$--domain may be a $((\overline{\star})_{_{\! f}},
\star_{_{\! f}})$--domain even if it is not integrally closed (cf.
Proposition \ref{finite-stable}).

\smallskip

\begin{cor} \label{cor:2.15} Let $D$ be a $\star$-domain. Assume that $D^\star$ is
flat over $D$.  Then $D$  is a $((\overline{\star})_{_{\! f}},
\star_{_{\! f}})$--domain.
\end{cor}
\begin{proof}
Let $\iota: D \hookrightarrow D^\star$ be the canonical embedding.
Then, $D^\star$ is a $\star_\iota$--domain, by Proposition
\ref{prop:2.7} (1). Since $D^\star$ is integrally closed
(Proposition \ref{ab}), we can apply Proposition
\ref{finite-stable} and get that $D^\star$ is a  ($(\overline{\
\!\star_\iota})_{_{\! f}}, (\star_\iota)_{_{\! f}})$--domain. By
using  also the flatness assumption of $D^\star$ over $D$, we have
$(E \cap F)^\star = ((E\cap F)D^\star)^\star = (ED^\star \cap
FD^\star)^\star = ((ED^\star)^\star \cap (FD^\star)^\star)^\star =
((ED^\star)^{\star_\iota} \cap
(FD^\star)^{\star_\iota})^{\star_\iota} = (ED^\star)^{\star_\iota}
\cap (FD^\star)^{\star_\iota} = (ED^\star)^{\star} \cap
(FD^\star)^{\star}= E^\star \cap F^\star$, for all $E,F \in
\ff(D)$.  The conclusion follows from Remark \ref{rk:2.10} (2).
 \end{proof}

\smallskip

 We conclude this section with a transfer-type result.

\smallskip

 \begin{prop} \label {iota-fin-stab} Let $D$ be an integral
domain and $\star$ a semistar operation on $D$. Let $\iota:D
\hookrightarrow D^\star$  be the canonical embedding. If $D$ is a
$((\overline{\star})_{_{\! f}}, \star_{_{\! f}})$--domain then
$D^\star$ is a $((\overline{\ \!\star_\iota})_{_{\! f}},
(\star_\iota)_{_{\! f}})$--domain.
\end{prop}
\begin{proof}  We prove the claim by using the equivalence (i)$\Leftrightarrow$(ii) of Remark \ref{rk:2.10} (2).
Let $E,F \in \ff(D^\star)$. There exist $E_0, F_0 \in \ff(D)$ such
that $E = E_0D^\star$ and $F = F_0 D^\star$. Then $(E \cap
F)^{\star_\iota} = (E_0D^\star \cap F_0D^\star)^{\star_\iota}  =
(E_0D^\star \cap F_0D^\star)^{\star}  \subseteq (E_0D^\star)^\star
\cap (F_0D^\star)^\star = E_0 ^\star \cap F_0 ^\star = (E_0 \cap
F_0)^\star \subseteq (E \cap F)^\star = (E \cap F)^{\star_\iota}$.
Thus $(E \cap F)^{\star_\iota} = (E_0D^\star)^\star \cap
(F_0D^\star)^\star = E^\star \cap F^\star
 =  E^{\star_\iota} \cap F^{\star_\iota} $.
\end{proof}
\medskip

\section{Pr\"ufer $\star$--multiplication domains and $\star$-coherence}

\smallskip

In this section, we look for conditions for a $\star$-domain to be
a P$\star$MD, by using coherent-like   conditions.

We say that a \it domain \rm $D$ is



%
%


%

\begin{enumerate}

\item[a)]  \;  \emph{$\star$--{extracoherent}}  if for all $E,F
\in \ff(D)$, with $0 \neq E\cap F$,  there exists $J \in \ff(D)$,
with $J \subseteq E\cap F$, such that $J^\star = E^\star \cap
F^\star$;

\item[b)]  \;   \emph{$\star$--coherent}  if for all
$E,F \in \ff(D)$, with $0 \neq E\cap F$,  there exists $J \in \ff(D)$,  such that $J^\star =E^\star \cap F^\star$ (i.e. $E^\star \cap F^\star$ is $\star$--finite \cite[page 650]{FP});

\item[c)]  \;   \emph{truly $\star$--coherent}  if for all
$E,F \in \ff(D)$, with $0 \neq E\cap F$,  there exists $J \in \ff(D)$,  such that $J^\star = (E \cap F)^\star$ (i.e. $E \cap F$ is $\star$--finite);

\item[d)]  \;    \emph{$\star$--quasi-coherent} if for each $F \in
\ff(D)$, $(D:F)^\star = G^\star$ for some $G \in \ff(D)$ (i.e. $(D:F)$ is $\star$--finite).

\end{enumerate}

\smallskip

\begin{rem}  \label{rk:3.1} Note that, without loss of generality, the properties a), b), c) and d) can be tested for all $E', F' \in \ff(D)$ and $E', F'$ ideals in $D$. As a matter of fact, if $E, F \in \ff(D)$, then for some nonzero elements $e, f \in D$, $eE, fF \subseteq D$, thus for $h := ef$ we have $E' :=hE, F' :=hF \in \ff(D)$ and $E',\ F' \subseteq D$.  Therefore if $J' \in \ff(D)$ is such that ${J'}^\star = {E'}^\star \cap {F'}^\star$ [respectively:
${J'}^\star= ({E'} \cap {F'})^\star$], then $(h^{-1}J')^\star= E^\star \cap F^\star$ [respectively: $(h^{-1}J')^\star=(E \cap F)^\star$].  Moreover, $J' \subseteq E' \cap F'$, then $J:= h^{-1}J' \subseteq h^{-1}E'  \cap h^{-1}F' = h^{-1}hE \cap h^{-1}hF= E\cap F$.
For d), for each $F \in \ff(D)$, let $f \in D$ be a nonzero element of $D$ such that $F' := fF \subseteq D$.
If $J' \in \ff(D)$ is such that ${J'}^\star = (D:F')^\star$, then it is easy to see that $J:= f^{-1}J' \in \ff(D)$ is such that $J^\star = (D:F)^\star$.
\end{rem}

\smallskip

Recall that given a semistar operation $\star$  on an integral domain $D$, $D$ is called a \emph{$\star$--Noetherian domain} if $D$ has the ascending chain condition on the quasi--$\star$--ideals (i.e. the nonzero ideals  $J$ of $D$  such that $J= J^\star \cap D$),  \cite[Section 3]{EFP}.

\smallskip

 \begin{exs} \label{examples3.2}    \bf(1) \rm
 \sl An integral domain $D$ and a
(semi)star operation $\star$ such that $D$ is
$\star$--quasi-coherent but it is neither  $\star$--coherent  nor
truly $\star$--coherent. \rm

  For $ \star = d$, the notions of
{$\star$--{extracoherent}} domain, { truly $\star$--coherent}
domain  and {$\star$--coherent} domain coincide with the classical
notion of coherent domain; the notion of {$d$--quasi-coherent}
domain coincides with the classical notion of quasi-coherent
domain \cite{BAD}.   Therefore it is sufficient to take a
quasi-coherent non-coherent domain (see \cite[Examples 4.4 and
5.3]{Gl00}). \rm

\bf (2) \sl  For $ \star =
v$, the notions of {$\star$--coherent} domain  and {$\star$--quasi-coherent}  domain coincide with the notion of $v$--coherent domain  \rm   \cite[Proposition 3.6]{FG}.

\bf (3) \sl  A $\star$--Noetherian domain (e.g. a Noetherian
domain) is  {truly $\star$--coherent}   (and {truly
$\star_{_{\! f}}$--coherent})   
\rm

Recall that in  a $\star$--Noetherian domain each nonzero
fractional ideal is  $\star_{_{\! f}}$--finite  \cite[Lemma 3.3]{EFP},
thus it is obvious that a $\star$--Noetherian domain is   truly
  {$\star$--coherent}    (or \ truly   {$\star_{_{\! f}}$--coherent}).


\bf (4) \rm \sl A Noetherian domain (thus, in particular, a {truly
$\star$--coherent})  is not ne\-cessarily a $\star$--{extracoherent}
domain. \rm (This fact  led
 us to use the terminology of    ``{extracoherent}''   for this type
of ``strong $\star$--coherence'', cf. also the following Theorem
\ref{coherent} (1).)

 In order to construct an example of the type announced above, we
start by recalling   that, for $\star = v$, even when  $D$ is
Noetherian,  $(E \cap F)^v $ maybe properly included in $ E^v \cap
F^v$. An explicit example was constructed in  \cite[page
4]{AACDMZ} as follows.   Let $K$ be a field and $X$ an
indeterminate,  set  $D:=K[\![X^3, X^4, X^5]\!]$, $E:=(X^3, X^4)$,
$F:=(X^3, X^5) $,  $ M:=(X^3, X^4, X^5)$. Note that $ E^v =
(D:(D:(X^3, X^4))) = (D: (X^{-3}D \cap X^{-4}D)) = (K[\![X^3, X^4,
X^5]\!] : K[\![X]\!]) =  (X^3, X^4, X^5) = M$;  similarly $F^v =
M$. Therefore  $(X^3) =(X^3)^v = (E \cap F)^v \subsetneq E^v \cap
F^v =  (X^3, X^4)^v  \cap (X^3, X^5)^v = M \cap M = M$.

\bf (5) \rm  Note that, even if a coherent domain is not
necessarily $\star$--{extracoherent} by (4),  e.g. for $\star =
v$,   however  it is an easy consequence of the definitions that
 \sl  a coherent domain (e.g.  a Pr\"ufer domain
\cite[Proposition 25.4 (1)]{gilmer}) is $\star$--{extracoherent}
domain,
  for each stable semistar operation $\star$. \rm

\end{exs}

\smallskip

\begin{lemma} \label{starf-coherent} Let $D$ be an integral domain and $\star$ a semistar operation on $D$. Then:
\begin{enumerate}
 \item[\bf (1) \rm]   The {$\star$--{extracoherent}} domains coincide with  the {$\star_{_{\! f}}$--{extracoherent}} domains and    the {$\star$--coherent} domains coincide with  the {$\star_{_{\! f}}$--coherent} domains.

  \item[\bf (2) \rm]  $D$ is a truly  {$\star_{_{\! f}}$--coherent} domain if and only if,  for all
$E ,F \in \ff(D)$, there exists $J \in \ff(D)$, with $J \subseteq
E \cap F$, such that $J^\star = (E \cap F)^\star$ (or,
equivalently, if  $E \cap F$ is  $\star_{_{\! f}}$--finite).  In
particular, a truly  {$\star_{_{\! f}}$--coherent} domain  is a truly
{$\star$--coherent} domain.

\item[\bf (3) \rm]  $D$ is a $\star_{_{\! f}}$--quasi-coherent domain if,
for each $F \in \ff(D)$, $(D:F)^\star = G^\star$ for some $G \in
\ff(D)$, with $G \subseteq (D:F)$ (or, equivalently, if $(D:F)$ is
$\star_{_{\! f}}$--finite).
 In particular, a {$\star_{_{\! f}}$--quasi-coherent} domain  is a {$\star$--quasi-coherent} domain.

\end{enumerate}
\end{lemma}
\begin{proof} (1) follows immediately from the definitions. (2) and (3) are straightforward consequences of \cite[Lemma 2.3]{FP}.
\end{proof}

\smallskip

\begin{rem} \label{quasi-coherent-f}
If $\star$ is a (semi)star operation on a domain $D$, then, in
particular, for each $F \in \ff(D)$, $(D:F)$ is a divisorial ideal
thus $(D:F) =(D:F)^\star = (D:F)^v$ \cite[Theorem 34.1 (3,
4)]{gilmer}, hence   $D$ is $\star_{_{\! f}}$--quasi-coherent if and only
if $D$ is $\star$--quasi-coherent (Lemma \ref{starf-coherent}
(3)).
 \end{rem}
\smallskip

\begin{prop} \label{sstarnew} Let  $\star$ be a semistar operation on an integral domain $D$.  Assume that $D$ is a $((\overline{\star})_{_{\! f}}, \star_{_{\! f}})$--domain (e.g. $\star$ is a stable semistar operation on $D$). Then:
\begin{enumerate}

\item[\bf (1) \rm]  If $D$ is  a $\star$--Noetherian domain  then
$D$ is $\star$--{extracoherent}.

\item[\bf (2) \rm] The notions of  truly $\star$--coherent domain and  $\star$--coherent domain coincide.

 \item[\bf (3) \rm]   $D$ is $\star$--quasi-coherent if and only if
$(D^\star:F)$ is $\star$--finite for each $F \in \ff(D)$.


  \item[\bf (4) \rm]  $\star$--coherent implies $\star$--quasi-coherent.
\end{enumerate}
\end{prop}
\begin{proof}
(1) Let $E,F \in \ff(D)$, with $0 \neq E\cap F$.
We have already observed that,  in a ($(\overline{\star})_{_{\! f}}, \star_{_{\! f}})$--domain,  $E^\star \cap F^\star = (E\cap F)^\star $.   Moreover, by the $\star$--Noetherianity, there exists $J \in\ff(D)$ such that $J^\star=(E \cap F)^\star $ and $J \subseteq E\cap F$ \cite[Lemma 3.3]{EFP}.

(2) and (3) \ are    obvious since in this situation,  for all $E,
F \in \ff(D)$,   $ (E\cap F)^\star = E^\star \cap F^\star$ (Remark
\ref{rk:2.10} (2));     similarly,  if $F= x_1D+x_2D+ \dots +
x_nD$,  in the present  situation we have $(D:F)^\star =
(x_1^{-1}D \cap x_2^{-1}D \cap  \dots \cap x_n^{-1}D)^\star  =
(x_1^{-1}D)^\star \cap (x_2^{-1}D)^\star  \cap \dots \cap
(x_n^{-1}D)^\star  = x_1^{-1}D^\star \cap x_2^{-1}D^\star  \cap
\dots \cap x_n^{-1}D^\star =(D^\star : F^\star) = (D^\star :F)$.

  (4)   If $F= x_1D+x_2D+ \dots + x_nD$, then $(D:F)^\star = (x_1^{-1}D \cap x_2^{-1}D \cap  \dots \cap x_n^{-1}D)^\star  = (x_1^{-1}D)^\star \cap (x_2^{-1}D)^\star  \cap  \dots \cap (x_n^{-1}D)^\star $ and $(x_1^{-1}D)^\star \cap (x_2^{-1}D)^\star  \cap  \dots \cap (x_n^{-1}D)^\star = G^\star$, for some $G \in \ff(D)$, by the $\star$--coherence of $D$.
\end{proof}

\begin{thm}    \label{coherent} Let $D$ be an integral domain and $\star, \star_1, \star_2$
semistar operations on $D$. Then:
\begin{enumerate}

\item[\bf (1) \rm]  If $D$ is $\star$--{extracoherent}  then $D$
is $\star$--coherent and truly $\star$--coherent.

\item[\bf (2) \rm] If  $D$ is  truly  $\star$--coherent  then $D$ is $\star$--quasi-coherent.





 \item[\bf (3) \rm]    Assume that $\star_1 \leq \star_2$. If
$D$ is truly
$\star_1$--coherent [respectively: $\star_1$--quasi-coherent] then $D$ is truly  $\star_2$--coherent, [respectively: $\star_2$--quasi-coherent].   \\
Assume, moreover, that $\star_2$ is stable.   If  $D$ is
$\star_1$--{extracoherent}  [respectively: $\star_1$--coherent]
then $D$ is $\star_2$--{extracoherent} [respectively:
$\star_2$--coherent].
\end{enumerate}
Let $\iota$ be the canonical embedding of $D$ in $D^{\star}$.
Then:
\begin{enumerate}

 \item[\bf (4) \rm]   $D$ is $\star$--coherent if and only if
$D^{\star}$ is $\star_\iota$--coherent.

\item[\bf (5) \rm]    If $D$ is
$\star$--{extracoherent}
 then $D^{\star}$ is
$\star_\iota$--{extracoherent}.

  \item[\bf (6) \rm]
Assume, moreover that $D$ is a  $((\overline{\star})_{_{\! f}},
\star_{_{\! f}})$--domain   (e.g. $\star$ is stable);   then $D$ is
truly $\star$--coherent [respectively: $\star$--quasi-coherent]
if and only if $D^{\star}$ is truly $\star_\iota$--coherent
[respectively: $\star_\iota$--quasi-coherent].

\end{enumerate}
\end{thm}
\begin{proof}
(1) follows from the definitions and from the fact that, in general,  $(E \cap F)^\star \subseteq E^\star \cap
F^\star$.

(2) Recall that,  if $F= x_1D+x_2D+ \dots + x_nD$, then $(D:F)^\star = (x_1^{-1}D \cap x_2^{-1}D \cap  \dots \cap x_n^{-1}D)^\star$, thus truly $\star$--coherent implies $\star$--quasi-coherent.

%



 (3)    In general, it is easy to see that if $\star_1 \leq
\star_2$ and  if an ideal is $\star_1$--finite, it is also
$\star_2$--finite. The second part of the statement follows from
the fact that if  $J \in \ff(D)$  is  such that $J^{\star_1} =
E^{\star_1} \cap F^{\star_1}$, then $J^{\star_2} =
(J^{\star_1})^{\star_2} = (E^{\star_1} \cap F^{\star_1})^{\star_2}
$. By the stability of $\star_2$, we have $(E^{\star_1} \cap
F^{\star_1})^{\star_2}  = (E^{\star_1})^{\star_2}  \cap
(F^{\star_1})^{\star_2} = E^{\star_2} \cap F^{\star_2} $.

  (4)   Assume that $D$ is $\star$--coherent. Let $E,F \in
\ff(D^\star)$ and let $E_0, F_0 \in \ff(D)$ be such that $E = E_0
D^\star$ and $F = F_0 D^\star$. Then, there exists $J_0 \in
\ff(D)$ such that $J_0^\star = E_0^\star \cap F_0^\star$.  Set
$J:=J_0D^\star \in \ff(D^\star)$, then $J^{\star_{\iota}}
=(J_0D^\star)^\star = J_0^\star = E_0^\star \cap F_0^\star = (E_0
D^\star)^\star \cap ( F_0D^\star)^\star = E^{\star_{\iota}}
\cap F^{\star_{\iota}}$.   Conversely, let $E,F \in \ff(D)$. Then
$ED^\star, FD^\star \in \ff(D^\star)$. It follows that there
exists $H \in \ff(D^\star)$ such that $H^\star = H^{\star_{\iota}}
=  (ED^\star)^{\star_{\iota}} \cap (FD^\star)^{\star_{\iota}}=
(ED^\star)^\star \cap (FD^\star)^\star$.  Let $H_0 \in \ff(D)$
such that $H = H_0D^\star$. Then $(H_0)^\star = (H_0D^\star)^\star
 = H^\star = (ED^\star)^\star \cap (FD^\star)^\star = E^\star
\cap F^\star$.

(5) Let $E,F,E_0,F_0,J_0$ like in the first part of the proof of
(4). Observe that, in this case, we can take  $J_0 \subseteq
E_0\cap F_0$. Then $J=J_0D^\star \subseteq (E_0\cap F_0) D^\star
\subseteq E_0D^\star \cap F_0D^\star = E \cap F$  and  we conclude
like in the first part of the proof of (4).

 (6)  The statement for the truly
coherent case follows from (4)  and Proposition \ref{sstarnew}
(2), since  $\star_\iota$ is a (semi)star operation on $D^\star$
and   if $D$ is a  $((\overline{\star})_{_{\! f}}, \star_{_{\!
f}}$)--domain, then $D^\star$ is a $((\overline{\
\!\star_\iota})_f, (\star_\iota)_f$)--domain  (Proposition
\ref{iota-fin-stab}).

\ \!  Now suppose that $D$ is   $\star$--quasi-coherent.  Let $ F
\in \ff(D^\star)$ and let $ F_0 \in \ff(D)$ be such that $F = F_0
D^\star$. We know that there exists $G_0 \in \ff(D)$ such that $
(D:F_0)^\star = G_0^\star$. Then by the assumption we have
$(D^\star : F_0D^\star)^\star = (D:F_0)^\star$. The conclusion is
now straightforward.


Conversely, assume that $D$ is  is $\star_\iota$--quasi-coherent.
Let $F_0  =x_1D+x_2D+\dots + x_nD   \in \ff(D)$. For $
F:=F_0D^\star \in \ff(D^\star)$, we know that there exists $G_0
\in \ff(D)$ such that  $G:=G_0D^\star \in \ff(D^\star)$ has the
property that $(D^\star: F_0D^\star)^{\star}=(D^\star:
F)^{\star_\iota} = G^{\star_\iota} = (G_0D^\star)^{\star_\iota}  =
G_0^\star$. Since  $D$ is a  $((\overline{\star})_{_{\! f}},
\star_{_{\! f}}$)--domain,
 we have that $(D^\star: F_0D^\star) =(D^\star: F_0) =
(D^{\star_{_{\! f}}}: (x_1D+x_2D+\dots + x_nD)) =
x_1^{-1}D^{\star_{_{\! f}}} \cap   x_2^{-1}D^{\star_{_{\! f}}}
\cap \dots \cap  x_n^{-1}D^{\star_{_{\! f}}} =  (x_1^{-1}D\cap
x_2^{-1}D\cap \dots \cap  x_n^{-1}D) ^{\star_{_{\! f}}} =
(D:F_0)^{\star_{_{\! f}}}$, therefore  $(D^\star:
F_0D^\star)^{\star} =  ((D:F_0)^{\star_{_{\! f}}})^\star =
(D:F_0)^\star$, thus we can conclude that $D$ is
$\star$--quasi-coherent.
\end{proof}

\begin{exs} \label {ex:3.7} \bf (1) \rm \sl A
$\star$--Noetherian domain (thus, in particular by Example
\ref{examples3.2}  (3), a truly $\star$--coherent domain) is not
necessarily $\star$--coherent. \rm

Let $D$ be a $2$-dimensional Noetherian domain. The integral
closure $D^\prime$ of $D$ is a $2$-dimensional Krull domain.
Clearly $D^\prime$ is not a Pr\"ufer domain (since a Pr\"ufer
Krull domain is a Dedekind domain \cite[Theorem 43.16]{gilmer} and
so $1$-dimensional). Thus, by \cite[Theorem 5.3.15]{glaz}, there
exists a proper overring $T$ of $D$ that is not coherent. Consider
the semistar operation $\star_{\{T\}}$ on $D$. Since $D$ is
Noetherian, it is obviously $\star_{\{T\}}$--Noetherian. Let
$\iota$ be the canonical embedding of $D$ in $T$. We have that $D$
is not $\star_{\{T\}}$--coherent, otherwise $T$ would be
$(\star_{\{T\}})_\iota$--coherent (Theorem \ref{coherent} (4)).
This is impossible, since  $(\star_{\{T\}})_\iota = d_T$, and
$d_T$--coherent  means coherent.

\bf (2) \rm \sl A  $v$-coherent domain  is not necessarily
$v$--extracoherent. \rm

Note that, since a $\star$--Noetherian domain is truly
$\star$--coherent (Example \ref{examples3.2} (3)), it is also
$\star$--quasi-coherent (Theorem \ref{coherent} (2)).
 Therefore, taking $\star = v$, a $v$--Noetherian domain (that is, a Mori
domain)
is $v$--quasi-coherent   or,  equivalently, $v$--coherent (Example
\ref{examples3.2} (2)). It follows that the Noetherian (in
particular, $v$--Noetherian) domain constructed in  Example
\ref{examples3.2} (4)  is  a $v$--coherent domain which is not
$v$--extracoherent   (and so, the notion of $\star$--coherence and
$\star$--extracoherence are distinct).
 \end{exs}


\begin{cor} \label{cor:3.7}
Let $D$ be an integral domain and $\star$ a semistar operation on
$D$. The following statements are equivalent:
\begin{enumerate}
\item[\bf (i) \rm] $D$ is $\til$--{extracoherent}.

\item[\bf (ii) \rm]  $D$ is truly $\til$--coherent.

 \item[ \bf (iii) \rm]  $D$ is
$\til$--coherent.
\end{enumerate}
\end{cor}
\begin{proof} We know already that (i)$\Rightarrow$(ii)$\Leftrightarrow$(iii) (Proposition \ref{sstarnew} (2) and Theorem \ref{coherent} (1)).

Assume that $D$ is truly $\til$--coherent. Since $\til$ is a semistar operation stable and of finite type, then for all $E,F \in \ff(D)$, there exists $J \in \ff(D)$ such that $J \subseteq E\cap F$ and $J^\til =(E\cap F)^\til = E^\til \cap F^\til$ \cite[Lemma 2.3]{FP}.
\end{proof}

\smallskip

 The next goal is to characterize the $\star$--extracoherence by using the other weaker $\star$--coherence-like conditions.  We start with an useful lemma

 \smallskip

\begin{lemma} \label{lemma:misc}
Let $D$ be an integral domain and $\star$ a semistar operation on
$D$. Assume that $D$ is $\star$--{extracoherent}. Then
 $D$ is a $ ((\overline{\star})_{_{\! f}},
\star_{_{\! f}})$--domain (or, equivalently,
 $D$ is a $(\widetilde{\star}, \star_{_{\! f}})$--domain).
\end{lemma}

\begin{proof}
Let   $E, F \in \ff(D)$.  Then there exists $J \in \ff(D)$, $J
\subseteq E \cap F$, such that $J^\star = E^\star \cap F^\star$.
Moreover,  obviously,   $J^\star \subseteq (E \cap F)^\star
\subseteq E^\star \cap F^\star$. Hence,  in particular,  $(E \cap
F)^\star = E^\star \cap F^\star$,  thus we conclude by Remark
\ref{rk:2.10} (2, (i)$\Leftrightarrow$(ii)).

 For the parenthetical statement, note that $\til \leq
(\overline{\star})_{_{\! f}} \leq \star_{_{\! f}}$, thus in the
present situation we need only to prove that $
((\overline{\star})_{_{\! f}}, \star_{_{\! f}})$--domain implies
$(\widetilde{\star}, \star_{_{\! f}})$--domain).   Since a
$\star$--{extracoherent} domain is also $\star_{_{\!
f}}$--{extracoherent}  (Lemma \ref{starf-coherent} (1)), we have
that in a $\star_{_{\! f}}$--{extracoherent} domain $
\widetilde{\star}  =(\overline{\ \!\star_{_{\! f}}})_{_{\! f}} =
 \star_{_{\! f}} $, by  what we have already proved  and
Remark \ref{rk:2.10} (3).
\end{proof}

\begin{rem} By using Lemma \ref{lemma:misc}, we can easily improve the result in Example
\ref{examples3.2}  (5)  and obtain that: \ \emph{A coherent domain
is $\star$--{extracoherent} if and only if it is a $
((\overline{\star})_{_{\! f}}, \star_{_{\! f}})$--domain.}   A
more precise statement will be proved in the following Proposition
\ref{prop:equiv-coherent}.
\end{rem}

\begin{prop} \label{prop:equiv-coherent}
Let $D$ be an integral domain and $\star$ a semistar operation on
$D$.
The following are equivalent:

\begin{enumerate}
\item[\bf (i) \rm] $D$ is $\star$--{extracoherent}.

 \item[\bf (i$_{\boldsymbol{f}}$) \rm]  $D$ is $\star_{_{\!
f}}$--{extracoherent}.

 \item[\bf ($\boldsymbol{\widetilde{\mbox{\bf i}}}$) \rm]  $D$ is
$\til$--extracoherent and a $(\widetilde{\star},\star_{_{\!
f}})$--domain.

  \item[ \bf (ii) \rm]  $D$ is truly $\star$--coherent and a
$(\widetilde{\star}, \star_{_{\! f}})$--domain.

 \item[\bf (ii$_{\boldsymbol{f}}$) \rm]  $D$ is truly
$\star_{_{\! f}}$--coherent and a $(\widetilde{\star},\star_{_{\!
f}})$--domain.

 \item[\bf ($\boldsymbol{\widetilde{\ \!\mbox{\bf ii}}\ \! }$)
\rm]   $D$ is truly $\widetilde{\star}$--coherent and a
$(\widetilde{\star},\star_{_{\! f}})$--domain.

 \item[\bf (iii) \rm]  $D$ is $\star$--coherent and a
$(\widetilde{\star},\star_{_{\! f}})$--domain.

 \item[\bf (iii$_{\boldsymbol{f}}$) \rm]  $D$ is $\star_{_{\!
f}}$--coherent and $(\widetilde{\star},\star_{_{\! f}})$--domain.

 \item[\bf ($\boldsymbol{\widetilde{\ \! \mbox{\bf iii}}\ \! }$)
\rm]   $D$ is $\tilde{\star}$--coherent and
$(\widetilde{\star},\star_{_{\! f}})$--domain.

\end{enumerate}
\end{prop}

\begin{proof}  The equivalences
(i)$\Leftrightarrow$(i$_f$) and (iii)$\Leftrightarrow$(iii$_f$) are in Lemma \ref{starf-coherent} (1).

    (ii$_f$)$\Leftrightarrow$($\widetilde{\mbox{ii}}$),
(iii$_f$)$\Leftrightarrow$($\widetilde{\mbox{iii}}$)  and
($\widetilde{\mbox{i}}$)$\Rightarrow$(i$_f$) are trivial.  Note
also that
($\widetilde{\mbox{i}}$)$\Leftrightarrow$($\widetilde{\mbox{ii}}$)$\Leftrightarrow$($\widetilde{\mbox{iii}}$)
by Corollary \ref{cor:3.7}.



 (ii)$\Leftrightarrow$(ii$_f$).  Observe that, when $\star_{_{\!
f}} = \widetilde{\star}$, for all  $E,F \in \ff(D)$,  we have $(E
\cap F)^{\star_{_{\! f}}} \subseteq (E \cap F)^{\star} \subseteq
E^{\star} \cap F^{\star} = E^{\star_{_{\! f}}} \cap F^{\star_{_{\!
f}}} = (E \cap F)^{\star_{_{\! f}}}$.
 Since in the present situation  $ (E \cap F)^{\star_{_{\! f}}}= (E \cap F)^{\star}$, it is clear that the notions of truly $\star$--coherent  and truly $\star_{_{\! f}}$--coherent  are equivalent.

(i$_f$)$\Leftrightarrow$(ii$_f$).
By the previous considerations, we already know that  (ii$_f$)($\Leftrightarrow(\widetilde{\mbox{ii}})\Leftrightarrow(\widetilde{\mbox{i}}))\Rightarrow$(i$_f$).

Conversely,  if $D$ is $\star_{_{\! f}}$--extracoherent, then
$D$ is truly $\star_{_{\! f}}$-coherent by Theorem
\ref{coherent} (1) and a $(\widetilde{\star},\star_{_{\!
f}})$--domain by Lemma \ref{lemma:misc} (and by
(i)$\Leftrightarrow$(i$_f$)).
\end{proof}

\medskip

The next goal is characterize the Pr\"ufer $\star$--multiplication
domains among the $\star$--domains using coherence-like
conditions.

\begin{prop} \label{p*md-coherent}
Let $D$ be an integral domain and $\star$ a semistar operation on
$D$. If $D$ is a P$\star$MD then $D$ is
$\tilde{\star}$--{extracoherent}.
\end{prop}
\begin{proof}
We claim that, in a P$\star$MD, for all $E,F \in
\boldsymbol{F}(D)$
$$((E+F)(E \cap
F))^{\tilde{\star}} = (EF)^\til.$$
Indeed, let $M \in \M(\star_{_{\! f}})$. Then $D_M$ is
a valuation domain and so $(E+F)(E \cap F) D_M = (ED_M + FD_M)(ED_M
\cap FD_M) = EFD_M$,  by  \cite[Theorem 25.2 (d) and Remark 25.3]{gilmer}. By the definition of $\widetilde{\star}$, we deduce the claim.

Now, if $E, F \in \ff(D)$, $EF$ is finitely generated, thus
$(EF(D:EF))D_M =D_M$, for each $M \in \M(\star_{_{\! f}})$ and so we
obtain that $EF$ is $\til$--invertible   \cite[Theorem 2.23]{FP}.

Therefore, $(E \cap F)$ is also $\til$--invertible
\cite[Lemma 2.1(2)]{FP} and, hence, $\til$--finite, i.e. $J^{\widetilde{\star}}= (E \cap F)^{\widetilde{\star}} = E^{\widetilde{\star}}  \cap F^{\widetilde{\star}}$, for some $J \in \ff(D)$, with $J \subseteq E \cap F$  \cite[Lemma 2.3 and  Proposition 2.6]{FP}.
\end{proof}

\begin{rem} \label{rk:p*md-coherent}
 Since in a P$\star$MD   it is known that $\widetilde{\star}= \star_{_{\! f}}$ \cite[Theorem 3.1 ((v)$\Leftrightarrow$(vi))]{FJS03b}, then  from the
previous Proposition \ref{p*md-coherent}  and Theorem
\ref{coherent} (1),  we deduce that a P$\star$MD is truly
$\widetilde{\star}$--coherent. Therefore, from  Proposition
\ref{prop:equiv-coherent},  we have that a P$\star$MD is also
$\star$--{extracoherent},  $\star_{_{\! f}}$--{extracoherent},
truly $\star$--coherent, truly $\star_{_{\! f}}$--coherent,
$\star$--coherent, $\star_{_{\! f}}$--coherent, and
$\widetilde{\star}$--coherent. Moreover, it is also
$\star$--quasi-coherent, $\star_{_{\! f}}$--quasi-coherent and
$\widetilde{\star}$--quasi-coherent, by Theorem \ref{coherent}
(2).

 Note   also that, from the fact that
 $\til \leq \overline{\star} \leq \star$, $\til \leq \star_{_{\! f}}$ and that $\overline{\star}$
is stable,  it follows that a P$\star$MD is
$\overline{\star}$--{extracoherent}  (Proposition
\ref{p*md-coherent} and Theorem \ref{coherent} (3))   and so truly
$\overline{\star}$--coherent, $\overline{\star}$--coherent and
$\overline{\star}$--quasi-coherent  (Theorem \ref{coherent} (1,
2)).

\end{rem}
\smallskip

 We are now in condition of proving the main theorem of this
section.

\smallskip
\begin{thm} \label{stardomain+coherent}
Let $D$ be an integral domain and $\star$ a semistar operation on
$D$. The following are
equivalent:

\begin{enumerate}
\item[\bf (i) \rm] $D$ is a P$\star$MD.

 \item[\bf (ii) \rm] $D$ is a $\star$--{extracoherent}
$\star$--domain.

 \item[\bf (ii$_{\boldsymbol{f}}$) \rm]  $D$ is a $\star_{_{\!
f}}$--{extracoherent} $\star$--domain.

 \item[\bf ($\boldsymbol{\widetilde{\ \!\mbox{\bf ii}}\ \!}$)
\rm]   $D$ is   a $\widetilde{\star}$--{extracoherent}
${\star}$--domain.

 \item[\bf (iii$_{\boldsymbol{f}}$) \rm]  $D$ is a truly
$\star_{_{\! f}}$--coherent $\star$--domain.

 \item[\bf ($\boldsymbol{\widetilde{\ \!\mbox{\bf iii}}\ \!}$)
\rm]   $D$ is   a  truly $\widetilde{\star}$--coherent
${\star}$--domain.





 \item[\bf ($\boldsymbol{\widetilde{\mbox{\ \!\bf iv}}\ \!}$)
\rm]   $D$ is   a $\widetilde{\star}$--coherent ${\star}$--domain.

\item[\bf (v$_{\boldsymbol{f}}$) \rm]   $D$ is  a $\star_{_{\!
f}}$--quasi-coherent ${\star}$--domain (or, equivalently, a
$\star_{_{\! f}}$--domain).

 \item[\bf ($\boldsymbol{\widetilde{\mbox{\bf v}}}$) \rm]   $D$
is a    $\til$--quasi-coherent ${\star}$--domain (or,
equivalently, a $\til$--domain).
\end{enumerate}
In particular, a quasi-coherent $\star$--domain is a  P$\star$MD.
\end{thm}
\begin{proof}
(i)$\Rightarrow$(${\widetilde{\mbox{ii}}}$)$\Leftrightarrow$(${\widetilde{\mbox{iii}}}$)$\Leftrightarrow$
(${\widetilde{\mbox{iv}}}$)$\Rightarrow$(${\widetilde{\mbox{v}}}$)$\Rightarrow$(v$_{f}$)
by Proposition \ref{basic} (3),   Proposition \ref{sstarnew} (2),
Corollary \ref{cor:3.7},  Proposition \ref{p*md-coherent} and
Theorem \ref{coherent} (1, 2, 3).

(i)$\Rightarrow$(ii)$\Leftrightarrow$(ii$_{f}$)  by Proposition
\ref{basic} (3), Remark \ref{rk:p*md-coherent} and Lemma
\ref{starf-coherent}  (1).

 (ii)$\Rightarrow$(iii),
(ii)$\Rightarrow$(iv)$\Leftrightarrow$(iv$_{f}$) and
(ii$_{f}$)$\Rightarrow$(iii$_{f}$)$\Rightarrow$(v$_{f}$) by
Theorem \ref{coherent} (1, 2) and Lemma \ref{starf-coherent} (1).

%

(v$_{f}$)$\Rightarrow$(i)  Let $F \in \ff(D)$. Then
$(FF^{-1})^\star = D^\star$, since $D$ is a $\star$--domain. By
the fact that $D$ is $\star_{_{\! f}}$--quasi-coherent,  we can
find $G \in \ff(D)$, with $G \subseteq (D:F)$,  such that
$G^{\star_{_{\! f}}}= (F^{-1})^{\star_{_{\! f}}}$ (Lemma
\ref{starf-coherent} (3)). Since $FG \in
 \ff(D)$ and $G \subseteq F^{-1}$, then we have
 $D^\star = (FF^{-1})^\star =  (F(F^{-1})^{\star_{_{\! f}}})^\star = (FG^{\star_{_{\! f}}})^\star=  (FG)^\star = (FG)^{\star_{_{\! f}}} \subseteq (FF^{-1})^{\star_{_{\! f}}} \subseteq D^\star$. Therefore  $F$ is
 $\star_{_{\! f}}$-invertible, and so $D$ is a P$\star$MD.

 For the parenthetical statements in  (v$_{f}$)  and (${\widetilde{\mbox{v}}}$)  see Proposition \ref{basic} (2) and Lemma \ref{starf-coherent} (3). The last claim follows from  Theorem \ref{coherent}  (3)   and  (v$_{f}$)$\Rightarrow$(i).
\end{proof}


\begin{rem} \label{rk:3.14}  \bf (1) \rm Note that, by  Remark \ref{quasi-coherent-f}, Theorem \ref{coherent} (1)  and the previous theorem,  \it when $\star$ is a (semi)star operation on $D$, then the following conditions are equivalent to each of the statements in Theorem \ref{stardomain+coherent}:

\begin{enumerate}
 \item[\bf (iii) \rm]  $D$ is a truly $\star$--coherent $\star$--domain.
 \item[\bf (iv) \rm]  $D$ is a $\star$--coherent (or, equivalently $\star_{_{\! f}}$--coherent, by Lemma \ref{starf-coherent} (1))
$\star$--domain.

 \item[\bf (v) \rm]  $D$ is a $\star$--quasi-coherent $\star$--domain.
 \end{enumerate} \rm

 \ec At the moment, we are unable to establish if the previous statements are also equivalent in the general semistar setting.  However, if $D$ is integrally closed or if $D^\star$ is flat over $D$ and if $(\overline{\star})_{_{\! f}}$ is stable, then (iii)$\Leftrightarrow$(iv) (and they are equivalent to each of the statements  of Theorem \ref{stardomain+coherent})  by Proposition \ref{finite-stable}, Corollary \ref{cor:2.15},  Lemma \ref{lemma:misc} and Proposition \ref{prop:equiv-coherent}.

 \bf (2) \rm From (1),  in particular,  we reobtain the
following result proved by Gabelli and Houston \cite[Proposition
3.2]{GH}:
$$
 D \mbox{ is a P$v$MD } \, \Leftrightarrow \,\,  D \mbox{ is a $v$--(quasi-)coherent $v$--domain}\,.
 $$
   Note that a similar characterization was given by Mott and Zafrullah  in \cite[Theorem 3.3]{MZ}:
 $ D$  is a P$v$MD if and only if   $D$  is essential and, for all nonzero $a, b \in D\, \,  aD \cap bD$ is finitely generated.  (Recall that an essential domain is a $v$--domain  by \cite[Proposition 44.13]{gilmer}
 and Remark \ref{rk:2.2} (5). On the other hand,  the condition called \it finite conductor \rm (for short, (FC)), i.e.  for all nonzero $a, b \in D\, \,  aD \cap bD$ is finitely generated is \sl  technically \rm weaker than the condition of quasi-coherence. However, in an essential domain, the condition (FC) is equivalent to the $v$--(quasi-)\-coherence \cite[Lemma 8]{Z78}.)

 \end{rem}

 \bigskip

Recall that a P$\star$MD which is  a $\star$--Noetherian domain is called a \emph{$\star$--Dedekind domain} \cite[Proposition 4.1]{EFP}.

\begin{cor}\label{cor:3.16}
Let $D$ be an integral domain and $\star$ a semistar operation on
$D$. The following are equivalent:

\begin{enumerate}
\item[\bf (i) \rm] $D$ is $\star$--Dedekind.

\item[\bf (ii)\rm ] $D$ is a $\star$--Noetherian $\star$--domain.
\end{enumerate}
\end{cor}

\begin{proof}
(i)$\Rightarrow$(ii) is obvious since a P$\star$MD is a
$\star$--domain (Proposition \ref{basic} (3)).

(ii)$\Rightarrow$(i) Since a $\star$--Noetherian domain is truly
$\star_{_{\! f}}$--coherent   (Example \ref{examples3.2} (3)), we
can apply Theorem \ref{stardomain+coherent}
((iii$_{f}$)$\Rightarrow$(i)).
\end{proof}


 \begin{exs}    \bf (1) \rm \sl A (semi)star operation $\star$
on an integral domain $D$ such that $D$ is
$\widetilde{\star}$--{extracoherent} but not $\star$--Noetherian
(and, so, not $\widetilde{\star}$--Noetherian \cite[Remark 3.6
(2)]{EFP}). \rm

Take a P$v$MD non-Krull (or, equivalently, non-Mori  \cite[Theorem
3.2]{K89}) domain $D$ (an explicit example is given next in
Example \ref{ex:4.1}). Then  $D$ is a $w$--{extracoherent}
$v$--domain (Theorem \ref{stardomain+coherent}), but $D$ is not
$v$--Noetherian, since a $v$--Noetherian  $v$--domain is a
$v$--Dedekind domain  (Corollary \ref{cor:3.16})  and
$v$--Dedekind coincides with Krull \cite[Remark 4.2 (1)]{EFP}.

\bf (2) \rm \sl A (semi)star operation $\star$ on an integral domain $D$ such that $D$ is not ${\star}$--quasi-coherent. \rm

 For $\star =v$,   take a $v$--domain $D$ which is not a P$v$MD
(Remark \ref{rk:2.2} (2)) then $D$ is not a $t$--quasi-coherent
(Theorem \ref{stardomain+coherent} ((v$_{f}$)$\Rightarrow$(i)))
or, equivalently, $D$ is not $v$--quasi-coherent   (Remark
\ref{quasi-coherent-f}).

 For $\star =d$, take any integrally closed non-P$v$MD $D$ and
apply \cite[Theorem 2]{Z78} to conclude that $D$ is not
quasi--coherent  (in fact, $D$ does not verify (FC)).

\end{exs}


\begin{rem}  \label{rk:3.8} \bf (1) \rm  In this section we have introduced a variety of coherence-like
de\-fi\-nitions for a semistar operation. But we were mainly interested in a ``strong form'' of semistar coherence,
 that we have called ``semistar--\-extracoherence'',  for obtaining a characterization of the P$\star$MDs
 in terms of coherence-like conditions (cf. Proposition \ref{p*md-coherent},  Theorem \ref{stardomain+coherent}).

We have also seen that all  the possibly different  coherence-like
notions introduced here coincide for a $\star$--domain, when
$\star$ is a (semi)star operation (Theorem
\ref{stardomain+coherent} and Remark  \ref{rk:3.14} (1)). However,
it seems to us that it would be interesting to investigate further
this subject in the general semistar setting and to study the
interconnections among the various coherence-like conditions in
some relevant situation (see, for instance, the following point
(3)).

 \bf (2) \rm  Note that we have introduced a notion of
``$\star$--coherence'', in order to be consistent
 with the definition of $v$--coherence already in the literature \cite{FG}, but we believe that
 the ``right'' definition of coherence in the semistar setting is what we called ``truly
 $\star$--coherence''. One of the reasons is that a $\star$--Noetherian domain is truly $\star$--coherent  but, in general,
 it is not $\star$--coherent  (Examples  \ref{examples3.2} (3) and \ref{ex:3.7} (1)).
 On the other hand, the fact that the notion of $v$--coherence works well in many situations is due to
  the fact that it  coincides with the $v$--quasi-coherence.

Recall that we have already shown that, in general,   a truly
$\star$--coherent domain is not a
  $\star$--coherent domain (Example \ref{ex:3.7} (1)). The following Example \ref{ex:3.18} shows conversely that  a $\star$--coherent domain is not a
   truly $\star$--coherent domain.

%
%

\bf (3) \rm  In this circle of ideas,  an open problem is related
to
 the specific cases of $v$--, $t$-- and $w$--operations.  Note that, by  Remark \ref{quasi-coherent-f}, Proposition \ref{starf-coherent} (1) and  \cite[Proposition 3.6]{FG},  we know already that: \ec
$$v\mbox{--coherent} \Leftrightarrow v\mbox{--quasi-coherent} \Leftrightarrow t\mbox{--quasi-coherent}
\Leftrightarrow t\mbox{--coherent}.$$

Therefore, by  Proposition \ref{sstarnew} (2) and Theorem \ref{coherent} (2), we have:
$$  w\mbox{--extracoherent}  \Leftrightarrow
\mbox{(truly) } w\mbox{--coherent}  \Rightarrow w\mbox{--quasi-coherent} \Rightarrow t\mbox{--quasi-coherent}
$$
and, by Theorem \ref{coherent}  (3),
$$
\mbox{(truly) } w\mbox{--coherent}  \Rightarrow \mbox{truly } t\mbox{--coherent} \Rightarrow \mbox{truly } v\mbox{--coherent} \Rightarrow v\mbox{--(quasi-)coherent}.$$

 \bf (Q-2) \sl  Is it possible  to
give examples for showing that the previous implications do not
invert~?
\end{rem}
\medskip

We end this section with the example announced in the previous remark.
\smallskip

 \begin{ex} \label{ex:3.18}  \sl A $\star$--coherent domain which is not truly $\star$--coherent (hence, not $\star$--extracoherent), for an a.b. semistar operation $ \star$ of finite type. \rm

 Let $k \subset K$ be a proper extension of fields and let $V$ a valuation domain of the form $K+M$, where $M$ is the nonzero maximal ideal of $V$
  such that $M=M^2$. (For instance, take $V$ to be a 1-dimensional nondiscrete valuation domain,  having value group $\mathbb R$.   Since $\mathbb R = 2 \mathbb R$, then clearly $M=M^2$. To produce examples of dimension greater than 1, take  $V$ having value group equal to the lexicographically ordered direct product $\mathbb R^n$, with $n\geq 2$.) \ec Set $D:= k+M$.

 Let $\star := \star_{\{V\}}$. Clearly $\star$ is an a.b. semistar operation of finite type on $D$. We claim that $D$ is not coherent but it is $\star$--coherent.

 Take $x \in K\setminus k$ and $m \in M$, $m \neq 0$, then we have $mD \cap mx D = mM$.  In order to prove this equality, we note that  $x$ is a unit in $V$  and $M$ is a common ideal in $D$ and $V$, thus obviously $xM =M$, and moreover
 $ mM= mM \cap mxM \subseteq mD \cap mxD$. On the other hand, if $y \in mD \cap mxD$, then
 $y = m(h_1 +m_1) = mx(h_2 +m_2)$,  with $h_1, h_2 \in k$ and $m_1, m_2 \in M$. Therefore  $ h_1- xh_2 = xm_2 -m_1 \in M$, and so $h_1 -xh_2 =0$. Since $x\in K\setminus k$, then necessarily $h_2 =0 =h_1$ and thus $y=mm_1 \in mM$.

 The equality $mD \cap xm D = mM$ and the fact that $M$ is not finitely generated (since $M=M^2 \neq 0$), implies that $D$ is not coherent (in fact, $D$ is not a finite conductor domain).  However, $(mD)^\star \cap (xm D)^\star =mV \cap xmV  = mV \cap mV = mV = (mD)^\star$. More generally if $E,F \in \ff(D)$ then $EV$ and $FV$ are principal and either $EV\cap FV =EV$ or $EV\cap FV =FV$, thus $D$ is clearly $\star$--coherent. But $D$ is not truly $\star$--coherent, since  $(mD \cap xm D)^\star = (mM)^\star = mM$, thus $FV = F^\star \neq  (mD \cap xm D)^\star = mM$, for each $F \in \ff(D)$.

 \end{ex}

\bigskip

\section{Pr\"ufer $\star$--multiplication domains and \texttt{H}$(\star)$--domains}

We recall from \cite[p. 651]{FP} that a domain $D$ with a semistar
operation $\star$ is an \it  \texttt{H}$(\star)$--domain \rm  if for each nonzero
integral ideal $I$ of $D$ such that $I^\star = D^\star$ there
exists a nonzero finitely generated ideal $J$, with $J \subseteq I$,  such that $J^\star = D^\star$ (i.e. $I$ is $\star_{_{\! f}}$--finite), in other words  $D$ is called an \texttt{H}$(\star)$--domain if $\F^\star = \F^{\star_{_{\! f}}}$.

When $\star =v$, the \texttt{H}$(v)$--domains coincide with the \texttt{H}--domains introduced by Glaz and Vasconcelos \cite[Remark 2.2 (c)]{GV77}.

It is obvious that each domain is an  \texttt{H}$(\star_{_{\! f}})$--domain, so the notion of  \texttt{H}$(\star)$--domain takes interest only when $\star$ is not of finite type.

%

In the next Proposition \ref{prop:2.9} we collect some
characterizations of the \texttt{H}$(\star)$--domains. Clearly a
$\star$--Noetherian domain is an \texttt{H}$(\star)$--domain
\cite[Lemma 3.3]{EFP}, thus we obtain in particular  that Mori
domains (e.g. Noetherian and Krull domains) are
\texttt{H}--domains. Houston and Zafrullah \cite[Proposition
2.4]{HZ88} proved, more generally,  that each $(t,v)$--domain (or
TV-domain in their terminology) is an \texttt{H}--domain.   Note
that a general class of \texttt{H}--domains which are not
$(t,v)$--domains was given  in \cite[Remark 2.5]{HZ88}.

\begin{ex} \label{ex:4.1}  \sl An  \texttt{H}$(\star)$--domain and P$\star$MD which is not a $\star$--Noetherian domain. \rm

 Clearly, for $\star = d$, a Pr\"ufer non-Dedekind domain provides an example of the type announced above.

A more elaborate example can be obtained by taking $\star = v$.
Let $K$ be a field and $X, Y$  indeterminates over $K$. Set $D:=
K[X] +YK(X)[Y]_{(Y)}$.  By the properties of the pullback
constructions  \cite{F}, $D$ is a 2-dimensional Pr\"ufer domain
with infinitely many maximal ideals, each of them invertible  (and
so divisorial) and with a unique height 1 prime ideal
$P:=YK(X)[Y]_{(Y)}$  which is also divisorial (since $P = (D:T)$,
where $T:=  K(X)[Y]_{(Y)}$)  and it is contained in all the
maximal ideals of $D$.  Clearly, in this case $\F^v = \F^t
=\{D\}$, thus $D$ is an  \texttt{H}$(v)$--domain. However,  $D$ is
not a Mori domain (in a Mori domain each nonzero element is
contained in finitely many maximal $t$--ideals  \cite[Proposition
2.2]{BG}, \cite[Th\'eor\`eme 4.2]{D})   and a Mori domain is
 precisely   a $v$--Noetherian domain  \cite[Section 3]{EFP} and
\cite[Theorem 2.1]{B}.  Note also that $D$ provides an explicit
example of an \texttt{H}$(v)$--domain which is not a
$(t,v)$--domain  \cite[Theorem 1.3 and Remark 2.5]{HZ88}.

\end{ex}

\begin{prop} \label{prop:2.9}  Let $\star$ be a semistar operation
on an integral domain $D$. The following conditions are
equivalent:
\begin{enumerate}
\item[\bf (i) \rm] $D$ is an  \texttt{H}$(\star)$--domain.


\item[\bf (ii) \rm] Each quasi--$\star_{_{\!f}}$--maximal ideal of $D$ is a
quasi--$\star$--ideal of $D$.

\item[\bf (iii) \rm]  For each $I \in\boldsymbol{{F}}(D)$,\ $I$ is $\star$--invertible
if and  only if $I$ is $\star_{\!_f}$--invertible.

\item[\bf (iii') \rm]  For each $I \in\boldsymbol{{F}}(D)$,\  if $I$ is $\star$--invertible
then  $I$ and $I^{-1}$ are $\star_{\!_f}$--finite.

\item[\bf (iv) \rm] $\mathcal{M}(\star_{\!_f}) = \mathcal{M}(\star)$.

\item[\bf (v) \rm]  $\mathcal{M}(\widetilde{\star}) = \mathcal{M}(\star)$.

\item[\bf (vi) \rm] The localizing system $\F^\star$ is finitely
  generated (i.e., $(\overline{\star})_{_{\! f}} = \overline{\star}$ \cite[Proposition 3.2]{FH2000}).

\item[\bf (vii) \rm]  $ \til =\overline{\star}$ (i.e. $D$ is a $(\til,\overline{\star} )$--domain).

\item[\bf (viii) \rm]  $\overline{\star} \leq \star_{_{\! f}}$.

\item[\bf (ix) \rm]   For each nonzero
prime ideal $P$ of $D$ such that $P^\star = D^\star$ there
exists a nonzero finitely generated  ideal $J$, with $J \subseteq P$,  such that $J^\star = D^\star$ (i.e. $P$ is $\star_{_{\! f}}$--finite).

\item[\bf (x) \rm]  $D$ is an  \texttt{H}$(\overline{\star})$--domain.

\end{enumerate}
\end{prop}

\begin{proof}
(i)$\Leftrightarrow$(ii) \cite[Lemma 2.7]{FP}.
(i)$\Leftrightarrow$(iii)$\Leftrightarrow$(iv)$\Leftrightarrow$(v)
\cite[Proposition 2.8]{FP}.

(i)$\Leftrightarrow$(vi) By definition,  the localizing system $\F^\star$ is the set of
the ideals $I$ of $D$ such that $I^\star = D^\star$. So, the definition of \texttt{H}$(\star)$--domain is equivalent to require that $\F^\star$ is finitely
generated.

(vi)$\Rightarrow$(vii)  Recall that  $\til = \star_{(\F^\star)_{_{\! f}}}$
(Proposition \ref{prop:loc1} (7)).  It is clear that, if
$\F^\star$ is finitely generated, then $\F^\star = (\F^\star)_{_{\! f}}$
\cite[Lemma 3.1]{FH2000}. Therefore  $\til = \star_{(\F^\star)_{_{\! f}}}
= \star_{\F^\star}= \overline{\star}$.

(vii)$\Rightarrow$(viii) It is clear, since $\til \leq \star_{_{\! f}}$.

(viii)$\Rightarrow$(iii) We have $\overline{\star} \leq \star_{_{\! f}} \leq
\star$, thus $ \F^\star \subseteq  \F^{\star_{_{\! f}}}  \subseteq  \F^{ \overline{\star}}$.   Since $ \F^\star =
 \F^{ \overline{\star}}$,  then obviously $ \F^\star= \F^{\star_{_{\! f}}} $. Therefore,  for an ideal $I \in \boldsymbol{{F}}(D)$,  we have  $II^{-1} \in  \F^\star$ if and only if  $II^{-1} \in  \F^{{\star}_{_{\! f}}}$; the conclusion now is straightforward.

 (vii)$\Leftrightarrow$(x)  Since $\overline{\overline{\star}} = {\overline{\star}} $ and $\widetilde{{\overline{\star}} } = \til$ \cite[Corollary 2.11 and Corollary 3.9]{FH2000}.

 (i)$\Rightarrow$(ix) is obvious.

 (ix)$\Rightarrow$(i) Assume that $D$ is not an \texttt{H}$(\star)$-domain, thus $\F^{\star_{_{\! f}}} \subsetneq
\F^\star $. It is easy to see that the nonempty set $\mathcal S :=
\F^\star \setminus \F^{\star_{_{\! f}}}$ is inductive.  Let $Q$ be a
maximal element in $\mathcal S $.  We claim that $Q$ is a prime
ideal of $D$. Suppose that $x,  y$ are two elements in $D\setminus
Q$ such that $xy \in Q$.  Then, by the maximality of $Q$ in
$\mathcal S$, we can find two finitely generated ideals of $D$,
$J'  \subseteq Q+xD$ and $J''  \subseteq   Q+yD$ such that ${J'}^{\star_{_{\! f}}}  =D^\star$ and $ {J''}^{\star_{_{\! f}}}
=D^\star$. On the other hand $J'J'' \subseteq Q^2 + xQ + yQ + xyD
\subseteq Q$, and $(J'J'')^{\star_{_{\! f}}} = ({J'}^{\star_{_{\! f}}}
{J''}^{\star_{_{\! f}}})^{\star_{_{\! f}}} = (D^\star)^\star = D^\star$, that
contradicts the fact that $Q \in \mathcal S$.

 Since $Q$ is a prime ideal and $Q\in \mathcal S \subset \F^\star$, then by assumption there exists a nonzero finitely generated ideal $J \subseteq Q$ such that $J^\star= Q^\star = D^\star$, i.e. $Q^{\star_{_{\! f}}} = D^\star$ or equivalently  $Q \in \F^{\star_{_{\! f}}} $, which is again  a contradiction.  \end{proof}

 \smallskip



\smallskip

\begin{cor} \label{star-dominio-H}
Let $D$ be an integral domain and $\star$ a semistar operation on
$D$. If $D$ is an \texttt{H}$(\star)$--domain the following conditions are
equivalent:

\begin{enumerate}

\item[\bf (i) \rm] $D$ is a $\star$-domain.

\item[\bf (ii)\rm ] $D$ is a P$\star$MD.

\end{enumerate}
\end{cor}

\begin{proof}
It is a straightforward consequence of Proposition
\ref{prop:2.9} ((i)$\Rightarrow$(iii)).
\end{proof}

Note that a P$\star$MD is not always an \texttt{H}$(\star)$--domain as the following example shows.

\begin{ex}
 \label{rem:H}  \sl A (semi)star operation $\star$ on an integral domain $D$ such that  $D$ is a P$\star$MD and a $((\overline{\star})_{_{\! f}}, \star_{_{\! f}})$--domain but not an \texttt{H}$(\star)$--domain and for which $\til  =\overline{\ \star_{_{\! f}}} = (\overline{ \star})_f= \star_{_{\! f}}  \lneq  \overline{\star} $. \rm

 Take a valuation domain $V$ with a non-divisorial  maximal ideal $M$ (e.g.  a rank 1 non-discrete valuation domain) and take  $\star = v$. Clearly $V$ is a P$v$MD, but not an
\texttt{H}$(v)$--domain, since the maximal ideal $M$ is a
$t$--ideal, but  not divisorial (Proposition \ref{prop:2.9} ((i)
$\Rightarrow$(iv))). Note that in this case $d =w\ (=
\widetilde{v}) = t\ (= v_f )$, thus $d=(\overline{v})_f = t
=\overline{v_f} $. Moreover, in a valuation domain, it is obvious
that every (semi)star operation is stable, thus in particular
$\overline{v}  =v$.  Finally,  by the previous considerations, it
follows that  $(d=) \ t  \lneq \overline{v} \ (=v)$,  since $M
\subsetneq M^v = V$.
\end{ex}

\smallskip

Example  \ref{rem:H} shows that the condition of being an
\texttt{H}$(\star)$--domain is too strong  to turn a $\star$--domain into a
P$\star$MD. On the other hand,  the condition of being an
\texttt{H}$(\star)$--domain is equivalent to the fact that  the subset
$\Inv(D, \star_{_{\! f}})$ of $\Inv(D, \star)$ coincides with  $\Inv(D, \star)$ (Proposition \ref{prop:2.9} ((i)$\Leftrightarrow$(iii))).
We can weaken condition (iii) of the previous  Proposition
\ref{prop:2.9}    and we
 call \it \texttt{I}$(\star)$--domain \rm an integral domain $D$ such that  $\Inv(D, \star) \cap \ff(D) = \Inv(D, \star_{_{\! f}}) \cap \ff(D)$.  Obviously  an \texttt{H}$(\star)$--domain is an \texttt{I}$(\star)$--domain and if $\star$ is a semistar operation of finite type on an integral domain $D$, then $D$ is  always  an \texttt{I}$(\star)$--domain.

%
%


\begin{prop} \label{weak}
Let $D$ be an integral domain and $\star$ a semistar operation on
$D$. The following conditions are equivalent:

\begin{enumerate}
\item[\bf (i) \rm] $D$ is an  \texttt{I}$(\star)$--domain.


\item[\bf (ii) \rm] If  $F\in \ff(D)$ is $\star$--invertible
 then $F^{-1}$ is $\star_{_{\! f}}$--finite.

 \item[\bf (iii) \rm] $D$ is an  \texttt{I}$(\overline{\star})$--domain.

\end{enumerate}
In particular, if $D$ is  $\star_{_{\! f}}$--quasi-coherent then $D$ is an  \texttt{I}$(\star)$--domain.
\end{prop}

\begin{proof}
(i)$\Rightarrow$(ii) If $F$ is $\star$--invertible, then it is
$\star_{_{\! f}}$--invertible. So, $F^{-1}$ is $\star_{_{\! f}}$--finite by
\cite[Proposition 2.6]{FP}.

(ii)$\Rightarrow$(i) It follows easily from \cite[Proposition
2.6]{FP}.

(iii)$\Leftrightarrow$(i).  Recall that a fractionary ideal $I$ is $\star$--invertible if and only if is $\overline{\star}$--invertible, since $\F^\star = \F^{\overline{\star}}$. By the equivalence (i)$\Leftrightarrow$(ii), we need to prove that if  $F\in \ff(D)$ is $\star$--invertible
 then $F^{-1}$ is $\star_{_{\! f}}$--finite if and only if  $F^{-1}$ is $({\overline{\star}})_f$--finite.

Let  $F^{-1}$ be $({\overline{\star}})_f$--finite.  If  $G \in \ff(D)$ is such that $G \subseteq (D:F)$ and $G^{{\overline{\star}}} = (D:F)^{{\overline{\star}}} $, then necessarily $G^{{{\star}}} = (D:F)^{{{\star}}} $, since ${\overline{\star}} \leq \star$.

 Let  $F^{-1}$ be $\star_{_{\! f}}$--finite.   If  $G \in \ff(D)$ is such that $G \subseteq (D:F)$ and $G^{{{\star}}} = (D:F)^{{{\star}}} $.  Since $F$ is $\star$--invertible then  $ (FG)^{{\star}} = (F(D:F))^{{{\star}}} = D^\star$, thus $FG \in \F^\star$. Since $\F^\star = \F^{\overline{\star}}$, then $ (FG)^{\overline{\star}}= D^{\overline{\star}} = (F(D:F))^{\overline{\star}}$.  Therefore $ ((FG)^{\overline{\star}} (D:F)^{\overline{\star}} )^{\overline{\star}} = (D^{\overline{\star}}(D:F)^{\overline{\star}})^{\overline{\star}}$, i.e. $ G^{\overline{\star}} =(D:F)^{\overline{\star}}$.

 Last statement is a straightforward consequence of  (ii)$\Rightarrow$(i) and of the fact that in a $\star_{_{\! f}}$--quasi-coherent,  for each $F \in \ff(D)$, $F^{-1}$ is $\star_{_{\! f}}$--finite.
\end{proof}

\medskip

\begin{rem}
\bf (1) \rm  Note that last statement of Proposition \ref{weak}
may not be improved by replacing  \texttt{I}$(\star)$--domain  with
 \texttt{H}$(\star)$--domain: Example \ref{rem:H}  provides  \sl a
  $t$--quasi-coherent domain \rm  (since P$v$MD, Theorem   \ref{stardomain+coherent}
((i)$\Rightarrow$(v$_{f}$)))  \sl
    which is not an \texttt{H}--domain. \rm

\bf (2) \rm The following Remark \ref{rk:4.10}  (2)  shows that the converse of the last statement of Proposition \ref{weak}  does not hold: it is easy to see that there exists \sl  an example of an \texttt{I}$(\star)$--domain which is not a $\star_{_{\! f}}$--quasi-coherent domain. \rm
\end{rem}

\medskip

The following result improves Corollary \ref{star-dominio-H} (cf.
also Theorem \ref{stardomain+coherent}
((i)$\Leftrightarrow$(v$_{f}$))  and last statement of Proposition
\ref{weak}).

\begin{cor} \label{star-dominio+weak}
Let $D$ be an integral domain and $\star$ a semistar operation on
$D$. The following conditions are
equivalent:

\begin{enumerate}

\item[\bf (i) \rm] $D$   is a $\star$-domain and an \texttt{I}$(\star)$--domain.
\item[\bf (ii)\rm ] $D$ is a P$\star$MD.
\item[\bf (iii)\rm ] $D$ is a P${\overline{\star}}$MD.

\end{enumerate}
\end{cor}

\begin{proof}  (i)$\Leftrightarrow$(ii)   Recall that, if  $F \in \ff(D)$, $I$ is $\star_{_{\! f}}$--invertible if and only if  $I$ is $\star$--invertible and  $F^{-1}$ is $\star_{_{\! f}}$--finite \cite[Proposition 2.6]{FP}. Therefore this equivalence  follows easily from  Proposition \ref{weak} ((ii)$\Leftrightarrow$(i)) and Proposition \ref{basic} (3).

(ii)$\Leftrightarrow$(iii) By Proposition \ref{weak} ((i)$\Leftrightarrow$(iii)) and Proposition \ref{basic} (4), this equivalence  is a straightforward consequence of  (ii)$\Leftrightarrow$(i).
\end{proof}

\begin{rem}  Note that Example \ref{rem:H}  provides also \sl  an example of an  \texttt{I}$(\star)$--domain
which is not an  \texttt{H}$(\star)$--domain, \rm  since, by Corollary \ref{star-dominio+weak}, a P$v$MD is also a
\texttt{I}$(v)$--domain.
\end{rem}

\smallskip

The considerations in Example \ref{rem:H}  and the equivalence
(i)$\Leftrightarrow$(iii$_f$)   in Theorem
\ref{stardomain+coherent} lead  us also to consider another weaker
condition of the property of being an
 \texttt{H}$(\star)$--domain (i.e. ($\til, \overline{\star})$--domain), namely  the notion of  $(\til, (\overline{\star})_{_{\! f}})$--domain.   It is obvious (like for the  \texttt{H}$(\star)$--domain case) that when $\star = \star_{_{\! f}}$ an integral domain is trivially a  $(\til, (\overline{\star})_{_{\! f}})$--domain.

 Recall that, it was shown in \cite{FH2000}
that $\til \leq (\overline{\star})_{_{\! f}}$ but in general these two semistar operations do not
coincide \cite[Example 3.11]{FH2000}.

\begin{prop} \label{H-coherent-tilde}
Let $D$ be an integral domain and $\star$ a semistar operation on
$D$. Assume that $D$ is either an \texttt{H}$(\star)$--domain or
truly $\star_{_{\! f}}$-coherent. Then $D$ is a $(\til, (\overline{\star})_{_{\! f}})$--domain.
\end{prop}

\begin{proof}

The case of an  \texttt{H}$(\star)$--domain is immediate since it is a $(\til,  \overline{\star})$--domain (Proposition \ref{prop:2.9} ((i)$\Rightarrow$(vii))).

So, assume that $D$ is truly $\star_{_{\! f}}$--coherent and let
$I$ be a nonzero finitely generated ideal of $D$. We have only to
show that $I^{\overline{\star}} \subseteq I^\til$. Let $x \in
I^{\overline{\star}}$. Then, there exists $J \in \F^{\star}$ such
that $xJ \subseteq I$. Since $J \subseteq x^{-1}I \cap D$, we have
that $x^{-1}I \cap D \in \F^\star$. Moreover, $x^{-1}I, D \in
\ff(D)$, so  by the assumption  $x^{-1}I \cap D$ is $\star_{_{\!
f}}$--finite. Let $H \subseteq x^{-1}I \cap D$ be a nonzero
finitely generated ideal such that $H^{\star_{_{\! f}}} = H^\star
= (x^{-1}I \cap D)^\star = D^\star$. It follows that $H \in
\F^{\star_{_{\! f}}}$. Moreover $xH \subseteq x(x^{-1}I \cap D)
\subseteq I$, and so $x \in I^\til$. Hence
$(\overline{\star})_{_{\! f}} = \til$.
\end{proof}

\begin{rem} \label{rk:4.10}  \bf (1) \rm   Note that
Example \ref{rem:H} shows that \sl  a $\til$--{extracoherent}
domain (or, a truly $\star_{_{\! f}}$--coherent domain) which is
also a $\star$--domain is  not necessarily an
\texttt{H}$(\star)$--domain. \rm

\bf (2) \rm  Note that, if $\star = \star_{_{\! f}}$, then
properties of being an \texttt{H}$(\star)$--domain, an
\texttt{I}$(\star)$--domain and a $(\til, (\overline{\star})_{_{\!
f}}) $--domain  are all trivially satisfied (recall that
$\overline{(\star_{_{\! f}})}= \til$  by Proposition
\ref{prop:loc1} (7)).  So, in particular, none of them implies the
$\star$--quasi-coherence and the $\star$--coherence (and so
neither the $\star$--ultracoherence nor the truly
$\star$--coherence). Indeed, it is enough to take \  $\star = d$ \
and consider
 an arbitrary  non-($d$--)quasi-coherent domain  (e.g. a non-finite conductor domain \cite{Gl00}).   This example shows, in particular, that  there exists \sl  an example of an \texttt{I}$(\star)$--domain which is not a $\star_{_{\! f}}$--quasi-coherent domain  (cf. also the following question \bf (Q-3) \rm  in Remark \ref{rem:4.14} (1)). \rm

  \end{rem}

\smallskip

Recall that we already know that:
$$
\begin{array}{rl}
\mbox{truly $\star_{_{\! f}}$--coherent domain } \Rightarrow & \hskip -5pt \mbox{$(\til, (\overline{\star})_{_{\! f}} )$--domain (Proposition \ref{H-coherent-tilde}),}\\
\mbox{$\star_{_{\! f}}$--quasi-coherent domain } \Rightarrow & \hskip -5pt  \mbox{\texttt{I}$(\star)$--domain (Proposition \ref{weak}), \:\:\:\:  and } \\
\mbox{truly $\star_{_{\! f}}$--coherent domain } \Rightarrow & \hskip -5pt  \mbox{$\star_{_{\! f}}$--quasi-coherent domain (Theorem \ref{coherent} (2)).}
\end{array}
$$
 The next goal is to relate the notions of $(\til,
(\overline{\star})_{_{\! f}} )$--domain and
\texttt{I}$(\star)$--domain.

\begin{prop} \label{tilde-weak}
Let $D$ be an integral domain and $\star$ a semistar operation on
$D$. If $D$ is a  $(\til, (\overline{\star})_{_{\! f}} )$--domain  then $D$ is an \texttt{I}$(\star)$--domain.
\end{prop}

\begin{proof}
Assume that $I$ is a finitely generated $\star$--invertible ideal. We already remarked that this is equivalent to the fact that $I$ is $\overline{\star}$--invertible. We want to show that $I$
is $\til$--invertible (and so, $\star_{_{\! f}}$--invertible \cite[Proposition 2.18]{FP}).   By \cite[Proposition 5.3 (2)]{FL06}, it is enough to show that $I$ is
$\til$--e.a.b.,  i.e., that if $(IF)^\til = (IG)^\til$ for some $F,G
\in \ff(D)$ then $F^\til = G^\til$.

Since, by assumption,  $\overline{\star}$ and $\til$ coincide on
$\ff(D)$, we have $(IF)^{\overline{\star}} =
(IG)^{\overline{\star}}$ and so $F^{\overline{\star}} =
G^{\overline{\star}}$, because $I$ is
$\overline{\star}$-invertible. Thus, since $F,G \in \ff(D)$, again
by the assumption, we have $F^\til = G^\til$.
\end{proof}


\begin{cor} \label{cor:4.10} Let $D$ be an integral domain and $\star$ a semistar operation on
$D$. Assume that $D$ is a $\star$--domain. The following are equivalent:
\begin{enumerate}
 \item[\bf (i$_{\boldsymbol{f}}$)]   $D$ is a truly $\star_{_{\! f}}$--coherent domain.
 \item[\bf (ii$_{\boldsymbol{f}}$)]   $D$ is a $\star_{_{\! f}}$--quasi-coherent domain.
\item[\bf (iii) \rm]  $D$ is a $(\til, (\overline{\star})_{_{\! f}} )$--domain.
\item[\bf (iv) \rm]  $D$ is an \texttt{I}$(\star)$--domain.
\end{enumerate}
\end{cor}
\begin{proof} We already observed that in general  (i$_f$)$\Rightarrow$(ii$_f$),  (iii)$\Rightarrow$(iv),  (i$_f$)$\Rightarrow$(iii) and  (ii$_f$)$\Rightarrow$(iv) (Theorem \ref{coherent} (2), Proposition \ref{weak}, Proposition \ref{H-coherent-tilde}   and Proposition \ref{tilde-weak}).

(iv)$ \Rightarrow$(i$_f$) When $D$ is a $\star$--domain, an
\texttt{I}$(\star)$--domain is a P$\star$MD and thus   (i$_f$)
holds (Corollary \ref {star-dominio+weak} ((i)$\Rightarrow$(ii))
 and Theorem \ref{stardomain+coherent}
((i)$\Rightarrow$(iii$_f$)).
\end{proof}

\smallskip

\begin{ex} \sl An  \texttt{I}$(\star)$--domain $D$ which is not a P$\star$MD (and so, in particular, which is not a $\star$--domain). \rm

Take any integral domain domain which is not Pr\"ufer and take $\star = d$.

\end{ex}

\smallskip

\begin{rem} \label{rem:4.14}  \bf (1) \rm The characterizations in Corollary \ref{star-dominio+weak}  lead to investigate more in depth the class of
 \texttt{I}$(\star)$--domains. In particular, in relation with Proposition \ref{weak},  Proposition \ref{tilde-weak} and Corollary \ref{cor:4.10},  it is natural to consider the following question (with $\star \neq d$ in order to avoid the trivial cases):

 \bf (Q-3) \sl
 Is it possible to find an example of an \texttt{I}$(\star)$--domain which is not a $(\til, (\overline{\star})_{_{\! f}} )$--domain or which is not a $\star_{_{\! f}}$--quasi-coherent domain (and so, also, not a $\star$--domain)~?  Is it possible to find such examples with $\star =v$ ? \rm

 \bf (2) \rm Since the notions of $\star$--domain and
$\overline{\star}$--domain  [respectively: the notions of
\texttt{I}$(\star)$--domain and
\texttt{I}$(\overline{\star})$--domain; the notions of  P$\star$MD
and  P$\overline{\star}$MD]  coincide (Proposition \ref{basic} (4)
[respectively: Proposition \ref{weak};  Remark \ref{rk:2.2} (4)]),  we have that, \it  if $D$ is a $\star$--domain, \it  the conditions (i)--(iv) of Corollary \ref{cor:4.10} are also  equivalent to each of the following:
\begin{enumerate}

 \item[\bf (\=i$_{\boldsymbol{f}}$) \rm]   $D$ is a truly $(\overline{\star})_{_{\! f}}$--coherent domain.
 \item[\bf (\=i\=i$_{\boldsymbol{f}}$) \rm]   $D$ is a  $(\overline{\star})_{_{\! f}}$--quasi-coherent domain.
\end{enumerate} \rm
 and, using also Theorem \ref{stardomain+coherent},
\it  these conditions are also equivalent to each of the following:

 \begin{enumerate}
\item[\bf (\~i) \rm]  \it $D$ is a (truly) $\til$--coherent domain.
\item[\bf ($\widetilde{\mbox{\bf ii}}$) \rm]   $D$ is a $\til$--quasi-coherent domain.
\end{enumerate}
and also to each of the following:
 \begin{enumerate}
 \item[\bf (j) \rm] $D$ is a $\star$--{extracoherent} (or, equivalently,  $\star_{_{\! f}}$--{extracoherent})
domain.
 \item[\bf (\=j) \rm] $D$ is a $\overline{\star}$--{extracoherent} (or, equivalently,  $(\overline{\star})_{_{\! f}}$--{extracoherent})
domain.

 \item[\bf ($\boldsymbol{\widetilde{\mbox{\bf j}}}$) \rm] $D$ is   a $\widetilde{\star}$--{extracoherent}
domain.



%

\end{enumerate} \rm
\end{rem}

\bigskip

 The equivalence (ii)$\Leftrightarrow$(iv) in the following
theorem provides  evidence for question \bf (Q-1) \rm in Remark
\ref{rk:2.10} (4).

\begin{thm} \label{star-domain-tilde}
Let $D$ be an integral domain and $\star$ a semistar operation on
$D$. The following are equivalent:

\begin{enumerate}

\item[\bf (i) \rm] $D$ is a $\star$--domain and a $(\til, {\star}_f )$--domain.

\item[\bf (\={i}) \rm] $D$ is a $\overline{\star}$--domain and a $(\til, (\overline{\star})_{_{\! f}} )$--domain.

\item[\bf  (ii) \rm] ${\star}$ is a.b. and  $D$ is a  $(\til, {\star}_f )$--domain.

\item[\bf (\={i}\={i}) \rm] $\overline{\star}$ is a.b. and $D$ is a  $(\til, (\overline{\star})_{_{\! f}} )$--domain.

\item[\bf  (iii) \rm] ${\star}$ is e.a.b. and  $D$ is a  $(\til, {\star}_f )$--domain.

\item[\bf (\={i}\={i}\={i}) \rm] $\overline{\star}$ is e.a.b. and $D$ is a  $(\til, (\overline{\star})_{_{\! f}} )$--domain.

\item[\bf (iv) \rm] $D$ is a P$\star$MD.

\end{enumerate}
\end{thm}
\begin{proof}
(i)$\Rightarrow$(ii)$\Rightarrow$(iii),  (\={i})$\Rightarrow$(\={i}\={i})$\Rightarrow$(\={i}\={i}\={i}) and (i)$\Rightarrow$(\={i}) follow from  Proposition \ref{ab},  Proposition \ref{basic} (4) and from the fact that $\til \leq (\overline{\star})_{_{\! f}} \leq \star_{_{\! f}}$ \cite[Proposition 3.6]{FH2000}.

(iii)$\Rightarrow$(iv) and (\={i}\={i}\={i})$\Rightarrow$(iv) Since $\star$ [respectively: $\overline{\star}$] is e.a.b. then ${\star}_f = \til$ [respectively:
$(\overline{\star})_{_{\! f}} = \til$] is (e.)a.b. \cite[Lemma 3.8]{FL06} and so $D$ is a P$\star$MD
by \cite[Theorem 3.1 ((v))$\Rightarrow$(i))]{FJS03b}.

(iv)$\Rightarrow$(i)  Recall that a P$\star$MD is an integral domain such that $\star_{_{\! f}} = \overline{(\star_{_{\! f}})} = \til$ and is (e.)a.b.  \cite[Theorem 3.1 ((i)$\Leftrightarrow$(vi))]{FJS03b}. Moreover we know that a P$\star$MD is a $\star$--domain (Proposition \ref{basic} (3)). \end{proof}


\begin{rem} \label{rk:4.14}  \bf (1) \rm Note that Example \ref{rem:H} provides \sl an example of a $(\til, {\star}_f )$--domain ( so, in particular, a $(\til, (\overline{\star})_{_{\! f}} )$--domain, hence also an \texttt{I}$(\star)$--domain, Proposition \ref{tilde-weak})    which is not a $(\til, \overline{\star})$--domain  (i.e. an \texttt{H}$(\star)$--domain),   \rm since a P$\star$MD is a $(\til, {\star}_f )$--domain  (Theorem \ref{star-domain-tilde} ((i)$\Leftrightarrow$(iv)).

 \bf (2) \rm Example \ref{b-domain} provides also \sl an example of an integral domain $D$ and a semistar operation $\star$ on $D$ such that $D$ is a $\star$--domain  but $D$ does not verify any of the  (equivalent)  conditions   of Corollary \ref{cor:4.10}  and Remark \ref{rem:4.14} (2).

\bf (3) \rm  Note that the implication (i)$\Rightarrow$(iv) in
Theorem \ref{star-domain-tilde}  could be also proved directly by
applying
 Proposition \ref{tilde-weak} and Corollary
\ref{star-dominio+weak} ((i)$\Rightarrow$(ii)).

 In Theorem \ref{star-domain-tilde}, after proving that  (i)$\Leftrightarrow$(ii)$\Leftrightarrow$(iii)$\Leftrightarrow$(iv), then the equivalences
(\=i)$\Leftrightarrow$(\=i\=i)$\Leftrightarrow$(\=i\=i\=i)$\Leftrightarrow$(iv) could be also obtained from
Corollary \ref{star-dominio+weak} ((ii)$\Leftrightarrow$(iii)) and from the fact that $\widetilde{\overline{\star}}= \til$ \cite[Corollary 3.9]{FJS03b}.

\bf (4) \rm
%
  It is easy to see that \it the statements in Theorem \ref{star-domain-tilde}  are also equivalent to  the following:
\begin{enumerate}

 \item[\bf (i$^\prime$) \rm] $D$ is a ${\overline{\star}}$--domain and a $(\til, {\star}_f )$--domain.

 \item[\bf (i${^{\prime\prime}}$) \rm] $D$ is a ${{\star}}$--domain and a $(\til, (\overline{\star})_{_{\! f}} )$--domain.

 \item[\bf (ii$^\prime$) \rm] ${\overline{\star}}$ is a.b. and $D$ is a  $(\til, {\star}_f )$--domain.

 \item[\bf (iii$^\prime$) \rm] ${\overline{\star}}$ is e.a.b. and $D$ is a  $(\til, {\star}_f )$--domain.
\end{enumerate} \rm

 But, in general, they are not equivalent to the following (weaker) statements:
\begin{enumerate}

 \item[\bf (ii$^{\prime\prime}$) \rm]  \it ${{\star}}$ is a.b. and $D$ is a  $(\til, ({\overline{\star}})_f )$--domain.  \rm

 \item[\bf (iii$^{\prime\prime}$) \rm] \it ${{\star}}$ is e.a.b. and $D$ is a  $(\til, ({\overline{\star}})_f )$--domain.
\end{enumerate} \rm

As a matter of fact, if $\star = \star_{_{\! f}}$, the condition
``$(\til, ({\overline{\star}})_f )$--domain''  is trivially
satisfied, but the condition ``${{\star}}$ is a.b.'' (or
``${{\star}}$ is e.a.b.'') does not imply P$\star$MD for lack of
stability.

\end{rem}

\bigskip



\begin{thebibliography}{10}

%
\bibitem{AACDMZ} D.~D. Anderson, D.~F. Anderson, D.~L. Costa,  D.~E. Dobbs, J. Mott and M. Zafrullah, \it Some characterizations of $v$--domains and related properties, \rm Colloquium Mathematicum \bf 58 \rm (1989), 1--9.
%
\bibitem{AC2000}
D.~D. Anderson and S.~J. Cook, \emph{Two star--operations and
their induced lattices}, Comm.  Algebra \textbf{28} (2000),
2461--2475.

\bibitem{AC2005}
D.~D. Anderson and Sharon M. Clarke, \emph{Star operations that distribute over finite intersections}, Comm.  Algebra \textbf{33} (2005),
2263--2274.
%
%
%
%
%

\bibitem{B} Valentina Barucci, \it Mori domains, \rm
in ``Non-Noetherian Commutative Ring Theory'' (Scott~T. Chapman and Sarah Glaz, eds.), Kluwer Academic Publishers, 2000, pp. 57--73.

\bibitem{BAD} Valentina Baurucci, David F. Anderson and David E. Dobbs, \it Coherent Mori domains and the principal ideal theorem, \rm Comm. Algebra \bf 15 \rm (1987), 1119-1156.


\bibitem{BG} Valentina Barucci and Stefania Gabelli, \it How far is a Mori domain from being a Krull domain?, \rm J. Pure Appl. Algebra \bf 45 \rm (1987), 101--112.


%
%
%
%
%


\bibitem{D}  Nicole Dessagnes, \it Intersections d'anneaux de Mori. Exemples, \rm C.R. Math. Rep. Acad. Sci. Canada \bf 7 \rm (1985), 355--360.


%



\bibitem{EFP}
Said ~El~Baghdadi, Marco Fontana, and Giampaolo Picozza.
\newblock \emph{Semistar {D}edekind domains.}
\newblock J. Pure Appl. Algebra {\bf 193} (2004), 27--60.

\bibitem{F} Marco Fontana, \it Topologically defined commutative rings, \rm Ann. Mat. Pura Appl. \bf 123 \rm (1980), 331--355.


\bibitem{FG}
Marco~Fontana and Stefania~Gabelli, \emph{On the class group and
the local class group of a pullback}, J. Algebra {\bf 181} (1996),
803--835.


%
\bibitem{FH2000}
Marco Fontana and James~A. Huckaba, \emph{Localizing systems and
semistar
  operations}, ``Non-Noetherian Commutative Ring Theory'' (Scott~T. Chapman and
  Sarah Glaz, eds.), Kluwer Academic Publishers, 2000, pp.~169--198.

\bibitem{FHP97} Marco Fontana, James~A. Huckaba, and Ira~J. Papick, \sl{Pr\"ufer domains}, \rm
 Marcel Dekker Inc., New York, 1997.
%
%
%
 \bibitem{FJS03b}{Marco Fontana, Pascual Jara and Eva Santos},
\emph{Pr{\"u}fer $\star$--multiplication domains and semistar
operations}.  J. Algebra Appl.  \bf 2 \rm (2003), 21--50.
%
%
%
 \bibitem {FL01a} Marco Fontana and K. Alan Loper, \it Kronecker function rings a general approach,
 \rm in  ``Ideal Theoretic Methods in Commutative Algebra" (D.D.
 Anderson and I.J. Papick, Eds.), M. Dekker Lecture Notes Pure Appl.
 Math.  \bf 220 \rm (2001), 189--205.

%
%
%
%
\bibitem{FL03} Marco Fontana and K.~Alan Loper, \emph{Nagata rings, {K}ronecker
function rings
  and related semistar operations}, Comm. Algebra \bf 31 \rm (2003), 4775--4805.


\bibitem{FL06}
Marco Fontana and K.~Alan Loper, \emph{A generalization of Kronecker function rings and Nagata rings}. Forum Math. (to appear).




\bibitem{FP}
Marco Fontana and Giampaolo Picozza. \emph{Semistar invertibility
on integral domains}, Algebra Colloq., {\bf 12} (2005), 645--664.




%

\bibitem{GH}
Stefania Gabelli and Evan Houston, \emph{Coherentlike conditions
in pullbacks},  Michigan Math. J. \textbf{44} (1997), 99--123.

%
%
%
%

\bibitem{gilmer}
Robert Gilmer, \sl{Multiplicative ideal theory},  \rm Marcel Dekker,
New York, 1972.

%
%
%
%

\bibitem{glaz}  Sarah Glaz,  \sl{Commutative coherent rings}, \rm Lecture Notes
in Mathematics, \bf 1371\rm, Springer-Verlag, Berlin, 1989.


\bibitem{Gl00}
Sarah Glaz, \it Finite conductor rings, \rm Proc. Am. Math. Soc. \bf 129 \rm (2000), 2833--2843.


\bibitem{GV77}
Sarah Glaz and Wolmer~V. Vasconcelos, \emph{Flat ideals. {II}},
Manuscripta
 Math. \textbf{22} (1977), 325--341.


 \bibitem{Gr} Malcom Griffin,\it Some results on $v$-multiplication rings, \rm Canad. J. Math. \bf 19 \rm (1967), 710--722.
%
%
%
  \bibitem{HK01} Franz Halter-Koch, \emph{Localizing systems, module
  systems and semistar operations}, J. Algebra \bf 238 \rm (2001),
 723--761.

%
%
%


\bibitem{H81} William Heinzer, \it An essential integral domain with a non-essential localization, \rm Canad. J. Math. \bf 33 \rm (1981), 400--403.

\bibitem{HO} William Heinzer and Jack Ohm, \it An essential ring which is not a $v$-multiplication ring, \rm Canad. J. math. \bf 25 \rm (1973), 856--861.



\bibitem{HMM} Evan Houston, Saroj Malik and Joe Mott, \it Characterizations of $\star$-multiplication domains, \rm Canad. Math. Bull. \bf 27 \rm (1984), 48--52.


\bibitem{HZ88} Evan Houston and Muhammad Zafrullah,
\emph{Integral domains in which each $t$--ideal is divisorial},
Michigan Math.  J. {\bf{35}} (1988), 291--300.



%
%
%
%
\rm
\bibitem{K89}
B.~G. Kang, \emph{The converse of a well-known fact about Krull domains}, J. Algebra \bf 124 \rm (1989), 284--289.

\bibitem{K89b} B.~G. Kang, \emph{Pr\"ufer {$v$}--multiplication domains and the
ring {$R[X]\sb
 {N\sb v}$}}, J. Algebra \textbf{123} (1989), 151--170.
%
%





  \bibitem{MO1993} Ry{\=u}ki  Matsuda and A. Okabe, {\it
On an AACDMZ question}, Math. J. Okayama Univ.  \bf 35 \rm (1993),
41--43.

\bibitem{matsumura}
Hideyuki Matsumura, \sl{Commutative ring theory}, \rm Cambridge
University Press,
  Cambridge, 1986.



\bibitem{MZ} Joe Mott and Muhammad Zafrullah, \emph{On Pr\"ufer $v$-multiplication domains}, \bf 35 \rm (1981), 1--26.

\bibitem{GP} Giampaolo Picozza,
\emph{Star operations on overrings and semistar operations}, Comm.
Algebra {\bf{33}} (2005), 2051--2073.

%
%
%
%
%
%
%
%
%
%
%




\bibitem{Uda}
Hirohumi Uda, \emph{L{CM}-stableness in ring extensions},
Hiroshima Math. J.
  \textbf{13} (1983),  357--377.


 \bibitem{WM97} Wang Fanggui and R. L. McCasland, \emph{On
 $w$--modules over strong Mori domains}, Comm.  Algebra  \bf 25 \rm (1997),
   1285--1306.
%
%
%
%
%
%
  \bibitem{Z78} Muhammad Zafrullah, \it On finite conductor domains, \rm Manuscripta Math. \bf 24 \rm (1978),  191--204.
%


\end{thebibliography}
\end{document}